\documentclass[seceqn,secthm]{elsart}

\usepackage{amsfonts}
\usepackage[dvips]{graphics,color}
\newcommand{\Z}{\mathbb{Z}}

\newcommand{\F}{\mathbb{F}}

\newcommand{\A}{\mathbb{A}}

\input amssym.def
\newsymbol\square 1003
\newsymbol\lozenge 1006
\newsymbol\surjarrow 1310
\newsymbol\injarrow 131A
\newsymbol\ltimes 226E
\newsymbol\rtimes 226F
\newsymbol\void 203F
\newsymbol\leftrightsquigarrow 1321
\def\END {$\square$}

\def\proclaim #1. #2\par #3\par {\medbreak
\noindent{\bf#1.\enspace}{\sl#2}\par\medbreak
\noindent{\bf Proof.} #3 \par}
\def\proclaimb #1. #2\par {\medbreak
\noindent{\bf#1.\enspace}{\sl#2}\par\medbreak}
 0
 2
 0
\newfont{\biggie}{cmbx10 scaled 1728}
\newtheorem{lemma}[equation]{Lemma}
\newtheorem{remark}[equation]{Remark}
\begin{document}
\begin{frontmatter}
\title{Cohomology of $Aut(F_n)$ in the $p$-rank two case}
\author{Craig A. Jensen}
\date{Submitted 11/99}

\smallskip

\address{Department of Mathematics, The Ohio State University\\
Columbus, OH 43210, USA\\ 
{\tt jensen@math.ohio-state.edu}}

\begin{abstract}

For odd primes $p$, we examine 
$\hat H^*(Aut(F_{2(p-1)}); \Z_{(p)})$,
the Farrell cohomology
 of the group of automorphisms of
a free group $F_{2(p-1)}$ on $2(p-1)$ generators, with coefficients
in the integers localized at the prime $(p) \subset \Z$.  This extends
results in \cite{[G-M]} by Glover and Mislin, whose calculations
yield
$\hat H^*(Aut(F_n); \Z_{(p)})$ for $n \in \{ p-1,p\}$
and is concurrent with work by Chen in \cite{[C]}
where he calculates
$\hat H^*(Aut(F_n); \Z_{(p)})$ for $n \in \{ p+1,p+2\}.$
The main tools used are
Ken Brown's ``normalizer spectral sequence'' from \cite{[B]},
a modification of Krstic and Vogtmann's proof of the contractibility
of fixed point sets for outer space in \cite{[K-V]}, and
a modification of the Degree Theorem of Hatcher and
Vogtmann in \cite{[H-V]}.
\end{abstract}


\end{frontmatter}

\section{Introduction}

Let $F_n$ denote the free group on $n$ letters and let
$Aut(F_n)$ and $Out(F_n)$ denote the automorphism group
and outer automorphism group, respectively, of $F_n$.
In \cite{[C-V]} Culler and Vogtmann defined a space on which
$Out(F_n)$ acts nicely called ``outer space''.  By
studying the action of $Out(F_n)$ on this space, various
people have been able to calculate the cohomology
of $Out(F_n)$ in specific cases.
In \cite{[Br]}, Tom Brady calculated the integral cohomology of 
$Out(F_3)$.
This remains even today the only complete (nontrivial) calculation of
the integral cohomology of $Out(F_n)$ or $Aut(F_n)$.
More recently, Hatcher in
\cite{[H]} and Hatcher and Vogtmann in \cite{[H-V]} have defined
a space on which $Aut(F_n)$ acts nicely called ``auter space''
and have used this to calculate the cohomology of
$Aut(F_n)$ in specific cases.

A ``Degree Theorem'' is introduced by Hatcher and Vogtmann
in \cite{[H-V]} which is a very useful tool for
simplifying cohomological calculations
concerning $Aut(F_n)$.
For example, they are able to derive linear stability
ranges for the integral cohomology of
$Aut(F_n)$ and able to calculate the rational cohomology of
$Aut(F_n)$ in some low dimensional cases.  In \cite{[JD]}, 
we modified this degree theorem and used it to calculate the 
cohomologies 
of some spaces having to do with $Aut(F_{2(p-1)})$.

Glover and Mislin \cite{[G-M]} calculated the
cohomology with coefficients in $\Z_{(p)}$ of $Out(F_n)$
for $n=p-1,p,p+1$.  In addition, Chen \cite{[Ch]} calculates
the integral cohomology of $Out(F_n)$ for $n=p+2$ and
$Aut(F_n)$ for $n=p+1,p+2$.  In each of the above cases,
the maximal $p$-subgroups of $Out(F_n)$ or $Aut(F_n)$ had $p$-rank 
one.
The case of $Aut(F_{2(p-1)})$ is the first one where
the maximal $p$-subgroups can have a higher $p$-rank, and it is
this case which we calculate here.

In this paper we use Brown's ``normalizer spectral sequence'' 
\cite{[B]},
to show 

\begin{figure}[here]
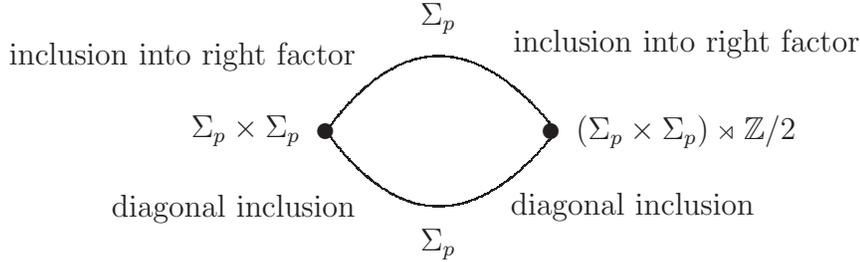

\caption{\label{gogps2} The graph of groups giving $\varsigma_p$}
\end{figure}

\begin{thm} \label{t9} Let $p$ be an odd prime, and $n=2(p-1)$. 
\newline
Then $\hat H^t(Aut(F_{n}); \Z_{(p)})$

$$\matrix{
\cong \left\{\matrix{
\Z/p^2 \oplus p(\Z/p) \hfill &t = 0 \hfill \cr
(p + [{3k \over 2}] - 1)\Z/p \hfill &|t| = kn \not = 0 \hfill \cr
\Z/p \hfill &t = 1 \hfill \cr
0 \hfill &t = kn + 1 > 1 \hfill \cr
[{3(k-1) \over 2}]\Z/p \hfill &t = -kn+1 < 0 \hfill \cr
H^{n-1}(\tilde Q_{p-1}; \Z/p) \oplus
[{3(k-1) \over 2}]\Z/p \hfill &t = kn-1 > 0 \hfill \cr
H^{n-1}(\tilde Q_{p-1}; \Z/p) \hfill &t = -kn-1 < 0 \hfill \cr
H^r(Q_{p}^\omega;\Z/p) \oplus
\sum_{i=0}^{p-1} H^r(\tilde Q_{i} \times Q_{p-1-i};\Z/p)
\hfill &t = kn + r, 2 \leq r \leq n-2 \hfill \cr
} \right. \hfill \cr}$$

Here the $(2k-2)$-dimensional space $Q_k$ is the quotient of 
the spine of
auter space $X_k$ by $Aut(F_k)$.  The $(2k-1)$-
and $(2p-4)$-dimensional spaces $\tilde Q_k$ and $Q_p^\omega$,
respectively, are the quotients of contractible spaces
(defined in this paper) $\tilde X_k$ and $X_p^\omega$ on 
which $F_k \rtimes Aut(F_k)$ and 
$(F_{p-2} \times F_{p-2}) \rtimes (\Z/2 \times Aut(F_{p-2}))$,
respectively, act properly with finite quotient.

Moreover, the submodule of $\hat H^t(Aut(F_{n}); \Z_{(p)})$
generated by all of the cohomology classes that are explicitly
listed above (that is, the ones not listed as coming from
quotient spaces) is a subring and is isomorphic as a ring to
$\hat H^t((\raisebox{-1.1ex}{\biggie *}\hbox{}_{i=1}^{p-1} 
\Sigma_p) \ast \varsigma_p; \Z_{(p)})$,
where
$\raisebox{-1.1ex}{\biggie *}\hbox{}_{i=1}^{p-1} 
\Sigma_p$ is the free product of $p-1$ copies of the 
symmetric group,  
$\varsigma_p$ is the fundamental group of the 
graph of groups pictured in Figure \ref{gogps2},
and 
$(\raisebox{-1.1ex}{\biggie *}\hbox{}_{i=1}^{p-1} 
\Sigma_p) \ast \varsigma_p$ is the free product
of these two groups.
All of the cohomology classes that are explicitly listed above,
with the exception of the one $\Z/p$ listed in dimension 1, 
are detected upon restriction to finite subgroups.  
Specifically, they are detected upon restriction
to the finite subgroups coming from stabilizers of marked 
graphs with underlying graphs (see Figure \ref{fig1})
$R_{2(p-1)}$, $R_k \vee \theta_{p-1} \vee R_{p-1-k}$
for $k \in \{1,\ldots,p-2\}$, $\Omega_{2(p-1)}$, and
$\Psi_{2(p-1)}$. 
\end{thm}

A quick note about our notation is appropriate here.
In general, groups without any additional
structure will be written using multiplicative notation
(e.g., $\Z/p \times \Z/p \cong (\Z/p)^2$) but modules like 
cohomology
groups will be written using additive notation
(e.g., $\Z/p \oplus \Z/p \cong 2(\Z/p)$.)  Hence
the case $t=0$ of
our main result above should be read as stating that
$$\hat H^0(Aut(F_{2(p-1)}); \Z_{(p)})
\cong \Z/p^2 \oplus (\Z/p \oplus \cdots \oplus \Z/p)$$
where there are $p$ copies of $\Z/p$ in
$(\Z/p \oplus \cdots \oplus \Z/p)$.

\hspace*{4.5cm} \input{cohom1.pic}

For the case of the prime $p=3$,
the cohomology groups of all of the
relevant quotient spaces are calculated in
Appendix \ref{appen} of this paper.  This allows us
to state concisely what the above result gives us in the case
$p=3$.  This calculation was also (independently)
done by Glover and Henn for both $Out(F_4)$
and $Aut(F_4)$.   

\begin{cor} \label{tapp}
\medskip

$$\hat H^t(Aut(F_{4}); \Z_{(3)}) = 
\left\{\matrix{
\Z/9 \oplus 3(\Z/3) \hfill &t = 0 \hfill \cr
([{3k \over 2}] + 2)\Z/3 \hfill &|t| = 4k \not = 0 \hfill \cr
\Z/3 \hfill &t = 1 \hfill \cr
0 \hfill &|t| \equiv 1,2  \hbox{ } (\hbox{mod } 4)
\hbox{ and } t \not = 1 \hfill \cr
[{3(k-1) \over 2}]\Z/3 \hfill &|t| = 4k-1 \hfill \cr
} \right.$$

Moreover, all of the cohomology classes, with the exception 
of the one $\Z/p$ listed in dimension 1, are detected upon 
restriction to finite subgroups.  Specifically, they are 
detected upon restriction to the finite subgroups coming 
from stabilizers of marked 
graphs corresponding to
$R_{4}$, $R_1 \vee \theta_2 \vee R_{1}$, $\Omega_{4}$, 
and $\Psi_{4}$.  Finally, there is a (ring) isomorphism
from the Farrell cohomology of $Aut(F_4)$ above to that
of $\Sigma_3 \ast \Sigma_3 \ast \varsigma_3$.
\end{cor}

In the next section, we introduce some useful spectral sequences,
in Section 3 we review the definitions of auter space
and its spine, and in Section 4, we compute the elementary abelian
$p$-subgroups of $Aut(F_n)$. For Section 5, we adapt
methods of Krstic to compute the normalizers of these
subgroups, in the section after that we
define certain contractible subcomplexes of
auter space that these normalizers act on,
and in the following section we use
the equivariant cohomology spectral sequence associated
to the action of the normalizers on these
contractible subcomplexes to compute
the cohomology of the normalizers.  Finally, Theorem \ref{t9} 
will be proved
in Section 8.  Appendix A states what Theorem \ref{t9} implies 
about the usual cohomology of $Aut(F_n)$, Appendix B 
briefly states some results about the cohomology of $Aut(F_l)$ 
in lower rank cases, and Appendix
C contains the proof of Corollary \ref{tapp}.

This paper is based on a dissertation written while the author 
was a student of Karen Vogtmann at Cornell, and the author would like
to thank Prof. Vogtmann for her help and advice.
The author would also like to thank Henry Glover for
his helpful comments on this paper.  Lastly, the author would like
to thank Hans-Werner Henn for pointing out errors in an earlier
version of this paper.

\section{Spectral sequences} \label{c1}

Let $G$ be a group acting cellularly on a
finite dimensional CW-complex $X$ such that the stabilizer
$stab_G(\delta)$ of every cell $\delta$ is finite and such
that the quotient of $X$ by $G$ is finite.  Further suppose that
for every cell $\delta$ of $X$, the group $stab_G(\delta)$ fixes
$\delta$ pointwise.  Let $M$ be a $G$-module. 
Recall (see \cite{[B]}) that the equivariant
cohomology groups of the $G$-complex $X$
with coefficients in $M$ are defined by
$$H_G^*(X; M) = H^*(G; C^*(X;M))$$
and that if in addition $X$ is contractible
(which will usually, but not always, be the case in this
paper) then
$$H_G^*(X; M) = H^*(G;M).$$

In \cite{[B]} Brown reviews the spectral sequence
for equivariant cohomology:
\begin{equation} \label{e1}
\tilde E_1^{r,s} = \prod_{[\delta] \in \Delta_n^r}
H^s(stab(\delta); M) \Rightarrow H_G^{r+s}(X; M)
\end{equation}
where $[\delta]$ ranges over the set $\Delta^r$ of orbits
of $r$-simplices $\delta$ in $X$.

If $M$ is $\Z/p$ or $\Z_{(p)}$ then
a nice property should be noted about the spectral sequence 
(\ref{e1}).
This property will greatly reduce the calculations we need to
go through, and in general will make concrete
computations possible.
Since each
group $stab(\delta)$ is finite, a standard restriction-transfer 
argument
in group cohomology yields that $|stab(\delta)|$ annihilates
$H^s(stab(\delta); M)$ for all $s > 0$.
(For examples of these sorts of
arguments see \cite{[A-M]} or \cite{[B]}.)
Since all primes not equal to
$p$ are divisible in $\Z/p$ or $\Z_{(p)}$,
this in turn shows that the $p$-part of
$|stab(\delta)|$ annihilates $H^s(stab(\delta); M)$ for
$s>0$.
In particular, if $p$ does not divide some $|stab(\delta)|$, then this
$[\delta]$ does not contribute anything to the spectral sequence
(\ref{e1}) except in the horizontal row $s=0$.  It follows that if our
coefficients are $\Z/p$ or $\Z_{(p)}$ then we are mainly just 
concerned
with the simplices $\delta$ which have ``$p$-symmetry''.

If $G$ is a group with finite virtual cohomological dimension (vcd)
and $M$ is a $G$-module, then Farrell cohomology groups
$$\hat H^*(G; M)$$
are defined in \cite{[F]}. 
For the basics about
Farrell cohomology, along with several useful properties,
see \cite{[B]} or \cite{[F]}. 

The equivariant cohomology spectral sequence for
Farrell cohomology is given by
\begin{equation} \label{e2}
E_1^{r,s} = \prod_{[\delta] \in \Delta^r}
\hat H^s(stab(\delta); M) \Rightarrow \hat H_G^{r+s}(X; M).
\end{equation}
and has analogous properties to those listed
for spectral sequence (\ref{e1}).

Ken Brown \cite{[B]} introduces another spectral sequence that can be
used to calculate $\hat H^*(G; \Z_{(p)})$.  It involves normalizers
of elementary abelian $p$-subgroups of $G$ and hence is often
called the ``normalizer spectral sequence.''  This is not the only
alternative available
to the standard spectral sequence for equivariant
cohomology, as Hans-Werner Henn \cite{[HN]} has created a
``centralizer spectral sequence'' that involves centralizers of
elementary abelian $p$-subgroups;
however, Brown's spectral sequence appears
to be the easiest to apply in our situation.

Let $p$ be a prime,
$G$ be a group with finite virtual cohomological dimension,
$\mathcal{A}$ be the poset of nontrivial elementary abelian
$p$-subgroups of $G$,
$\mathcal{B}$ be the poset of conjugacy classes of
nontrivial elementary
abelian $p$-subgroups of $G$, and $|\mathcal{B}|_r$ be the
set of $r$-simplices in the realization $|\mathcal{B}|$.
Brown's normalizer spectral sequence is
\begin{equation} \label{er4}
E_1^{r,s} = \prod_{(A_0 \subset \cdots \subset A_r) \in 
|\mathcal{B}|_r}
\hat H^s( \bigcap_{i=0}^r N_{G}(A_i); \Z_{(p)})
\Rightarrow \hat H^{r+s}(G; \Z_{(p)})
\end{equation}

We will use the above normalizer spectral sequence to calculate the
cohomology of $Aut(F_n)$.  
The cohomology groups of the normalizers, which are required as input
into the normalizer spectral sequence, will be computed using 
equivariant cohomology spectral sequences via their actions on certain
fixed point subspaces of auter space.

\section{$Aut(F_n)$ and auter space}

We review some basic properties and definitions of 
the automorphism group $Aut(F_n)$ of a free group $F_n$ of rank
$n$ (where $n=2(p-1)$ for our work.)   
Most of these can be found in \cite{[C-V]},
\cite{[H-V]}, \cite{[S-V1]}, \cite{[S-V2]}, and \cite{[Z]}.
Let $(R_n,v_0)$ be the
$n$-leafed rose,
a wedge of $n$ circles. We say a basepointed graph $(G,x_0)$
is \emph{admissible} if it has no free edges, all vertices except the
basepoint
have valence at least three, and there is a basepoint-preserving 
continuous
map $\phi\colon R_n \to G$ which induces an isomorphism on $\pi_1$. 
The
triple $(\phi,G,x_0)$ is called a \emph{marked graph}.  Two marked 
graphs
$(\phi_i,G_i,x_i) \hbox{ for } i=0,1$  are \emph{equivalent}  if 
there is
a homeomorphism $\alpha \colon (G_0,x_0) \to (G_1,x_1)$ such that
$ (\alpha\circ\phi_0)_\# = (\phi_1)_\# : \pi_1(R_n,v_0) \to 
\pi_1(G_1,x_1)$.
Define a partial order on the set of all equivalence classes of
marked graphs by setting $(\phi_0,G_0,x_0) \leq
(\phi_1,G_1,x_1)$ if $G_1$ contains a \emph{forest}
(a disjoint union of trees in $G_1$ which contains all of the vertices
of $G_1$) such that collapsing each tree in
the forest to a point yields $G_0$, where the collapse is compatible 
with
the maps $\phi_0$ and $\phi_1$.

From \cite{[H]} and \cite{[H-V]} we have that $Aut(F_n)$ acts
with finite stabilizers on
a contractible space $X_n$, called the spine of auter space.
The space $X_n$ is the geometric realization of the poset of
marked graphs that we defined above.
Let $Q_n$ be the quotient of $X_n$ by
$Aut(F_n)$.
Note that the CW-complex $Q_n$ is not necessarily a simplicial
complex.
Since $Aut(F_n)$ has a torsion free subgroup of finite index
\cite{[H]} and it acts on the contractible,
finite dimensional space
$X_n$ with finite stabilizers and finite quotient,
$Aut(F_n)$ has finite vcd.  Thus it makes sense to
talk about its Farrell cohomology, and to apply
the normalizer spectral
sequences  calculate its cohomology.

For $n=0$, we set
$X_0 = Q_0 = \{\hbox{a point}\}$
as a notational convenience.

We will use spectral sequence (\ref{er4}) to calculate the 
cohomology of $Aut(F_n)$:
\begin{equation} \label{e3}
\matrix{E_1^{r,s} = \prod_{(A_0 \subset \cdots \subset A_r) \in 
|\mathcal{B}|_r}
\hat H^s( \bigcap_{i=0}^r N_{Aut(F_n)}(A_i); \Z_{(p)}) 
\hspace*{3cm}\hfill \cr
\hfill \Rightarrow \hat H^{r+s}(Aut(F_n); \Z_{(p)}) \cr}
\end{equation}
In order to use the above spectral sequence, we must 
classify the elementary abelian $p$-subgroups of $Aut(F_n)$.
From Zimmerman's realization theorem \cite{[Z]} 
(cf. Culler \cite{[C]}), any finite
subgroup of $Aut(F_n)$ is realized on a marked graph, 
where a subgroup $G$ is said to be {\em realized} by
a specific marked graph $\eta : R_n \to \Gamma$
if it is contained in the stabilizer of that marked graph.
Smillie and Vogtmann \cite{[S-V1]}
examined the structure of these stabilizers in
detail, and we list their results here.
Consider a given $r$-simplex
$$(\phi_r,\Gamma_r,x_r) > \cdots  > (\phi_1,\Gamma_1,x_1)
> (\phi_0,\Gamma_0,x_0)$$
with corresponding forest collapses
$$(H_r \subseteq \Gamma_r), \ldots, (H_2 \subseteq \Gamma_2), 
(H_1 \subseteq \Gamma_1).$$
For each $i \in {0,1, \ldots, r}$, let $F_i$ be the inverse image
under the map
$$\Gamma_r \to \cdots \to \Gamma_{i+1} \to \Gamma_i$$ of forest 
collapses,
of the
forest $H_i$.  That is, we
have
$$F_r \subseteq \cdots \subseteq F_2 \subseteq F_1 \subseteq 
\Gamma_r.$$
The stabilizer of the simplex under
consideration is isomorphic to the group
$Aut(\Gamma_r,F_1,\ldots,F_r,x_r)$ of basepointed automorphisms of 
the graph
$\Gamma_r$ that respect each of the forests $F_i$.  For example,
the stabilizer of a point $(\phi,\Gamma,x_0)$ in $X_n$ is isomorphic 
to
$Aut(\Gamma,x_0)$.  If a marked graph $(\phi,\Gamma,x_0)$
realizes a subgroup $G$, we can
think of $\Gamma$ as a graph with an action of $G$ on it.
Hence a first step toward calculating the elementary
abelian $p$-subgroups of $Aut(F_n)$ will be finding out which
graphs have $p$-symmetry.

\section{Conjugacy classes of $p$-subgroups of $Aut(F_n)$} \label{c2}

From \cite{[B-T]} and \cite{[M]},
we see that $p^2$ is an upper bound for the order of any $p$-subgroup
of $Aut(F_n)$.  (Recall that we set $n=2(p-1)$ earlier.)
Every nontrivial
$p$-subgroup of $Aut(F_n)$ is isomorphic to either $\Z/p$ or
$\Z/p \times \Z/p$
(see, for example, Smillie and Vogtmann in \cite{[S-V2]}.)
 
We want to find all conjugacy classes of (elementary abelian)
$p$-subgroups of $Aut(F_n)$.
By Zimmerman's realization theorem \cite{[Z]}, 
we can do this by analyzing marked
graphs with $p$-symmetry, so that the stabilizers of these marked
graphs have elements of order $p$. 
The action of $Aut(F_n)$ is transitive on marked graphs with the
same underlying graph.
Since we are only interested in
conjugacy classes of subgroups, we will just examine the underlying
graphs. 

We first calculate which
graphs $\Gamma$ with a $\Z/p$
action on them have $\pi_1(\Gamma) \cong F_n$.
To simplify our calculations, we will consider only 
reduced $\Z/p$-graphs, where a
$\Z/p$-graph
$\Gamma$ is \emph{reduced} if it contains no $\Z/p$-invariant
subforests.

\hspace*{1.5cm}\input{cohom2.pic}
\begin{figure}[here]
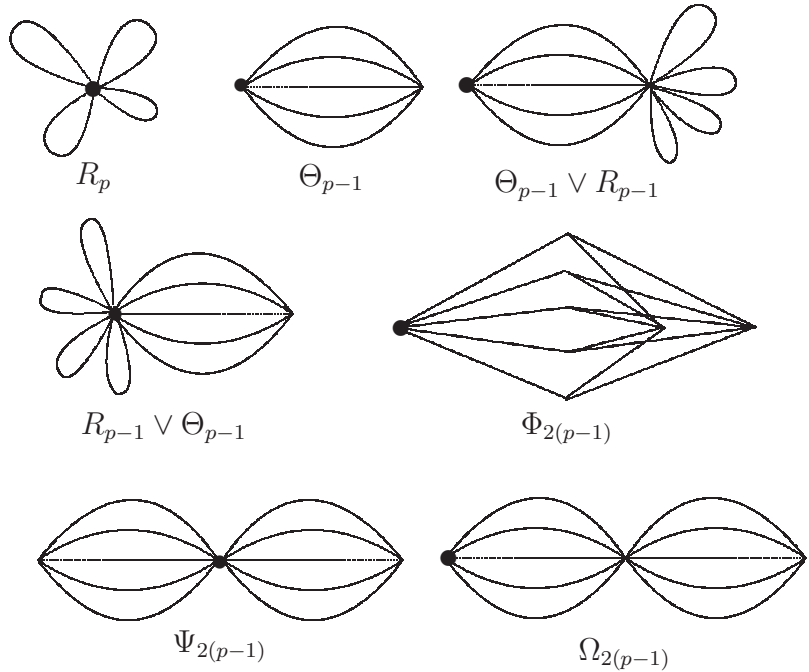

\caption{\label{fig1} Some graphs with $p$-symmetry}
\end{figure}

Some preliminary definitions of a few common graphs are in order.
Let $\Theta_{p-1}$ be the graph with two vertices and $p$ edges, 
each of
which goes from one vertex to the other (see Figure \ref{fig1}.)
Say the ``leftmost vertex''
of $\Theta_{p-1}$ is the basepoint.  Hence when we write
$\Theta_{p-1} \vee R_{p-1}$ then we are stipulating that the rose
$R_{p-1}$ is attached to the non-basepointed vertex of $\Theta_{p-1}$,
while when we write $R_{p-1} \vee \Theta_{p-1}$ then we are saying
that the rose is attached to the basepoint of $\Theta_{p-1}$.
Let $\Phi_{2(p-1)}$ be a graph with $3p$ edges $a_1, \ldots, a_p,$
$b_1, \ldots, b_p,$ $c_1, \ldots, c_p,$ and $p+3$ vertices
$v_1, \ldots, v_p, x, y, z.$  The basepoint is $x$ and each of the
edges $a_i$ begin at $x$ and end at $v_i$.  The edges $b_i$ and $c_i$
begin at $y$ and $z$, respectively, and end at $v_i$.  Note that there
are obvious actions of $\Z/p$ on $\Theta_{p-1}$ and $\Phi_{2(p-1)}$,
given by rotation, and that these actions are unique up to conjugacy.
Let $\Psi_{2(p-1)}$ be the graph obtained from $\Phi_{2(p-1)}$ by
collapsing all of the edges $a_i$ to a point.  Let $\Omega_{2(p-1)}$
be the graph obtained from $\Phi_{2(p-1)}$ by collapsing either the
edges $b_i$ or the edges $c_i$ (the resulting graphs are isomorphic)
to a point.  Note that the only difference between $\Psi_{2(p-1)}$ and
$\Omega_{2(p-1)}$ is where the basepoint is located.

\begin{lemma} \label{tr1}
Let $p$ be an odd prime and $n=2(p-1)$,
and let $\Gamma$
be a reduced (basepointed) graph
with a
nontrivial $\Z/p$-action, where $\pi_1(\Gamma) \cong F_n$.
Let $e$ be an edge of $\Gamma$ which is moved
by the $\Z/p$ action.  
Then the orbit $e_1, e_2, \ldots, e_p$ of
the $1$-cell $e$ under the $\Z/p$-action on the CW-complex
$\Gamma$
forms either a rose $R_p$ or a $\Theta$-graph $\Theta_{p-1}$.
\end{lemma}

\begin{pf*}{Proof.} 
Since $\Gamma$ is reduced,
the edges $e_i$ have more than one endpoint in common,
else they form an invariant subforest. They also cannot form a 
$p$-gon, else
we could take a minimal path from $e_1$ to the basepoint, consider
its orbit under $\Z/p$, and find in those edges a star of $p$ edges
that could be collapsed. \END \end{pf*}

\begin{prop} \label{t1}
Let $p$ be an odd prime and $n=2(p-1)$.
The only reduced (basepointed) graphs $\Gamma$
with a
nontrivial $\Z/p$-action
and $\pi_1(\Gamma) \cong F_n$ are $R_n$,
$R_k \vee \Theta_{p-1} \vee R_{p-1-k}$, $\Psi_n$, and $\Omega_n$.
\end{prop}

\begin{pf*}{Proof.} First, suppose that $\Gamma$ has only one 
nontrivial $\Z/p$-orbit of
edges $$e_1, e_2, \ldots, e_p.$$
From the lemma, 
$\Gamma$ is $R_{2(p-1)}$ with $\Z/p$ rotating $p$ of the leaves
and leaving the other $p-2$ fixed, or $\Gamma$ is
of the form $R_k \vee \Theta_{p-1} \vee R_{p-1-k}$
where $\Z/p$ rotates the edges of the $\Theta$-graphs and leaves
the roses at either end of the $\Theta$-graph fixed.

Second, suppose $\Gamma$ has more than one
nontrivial $\Z/p$-orbit of unoriented edges.  Take two
distinct orbits $e_1, e_2, \ldots, e_p$ and
$f_1, f_2, \ldots, f_p$.  As in the previous paragraph,
we can use the fact that $\Gamma$ is reduced
and Lemma \ref{tr1} to get that the $e_i$
and $f_i$ either form roses or $\Theta$-graphs.  If either one
of them is a rose, then the rank of $\pi_1(\Gamma)$ is at
least $2p-1$, which is a contradiction.  So both form
$\Theta$-graphs.  Since $\pi_1(\Gamma) \cong F_{2(p-1)}$, 
it follows that $\Gamma$ is either $\Psi_n$ or $\Omega_n$. \END 
\end{pf*}

We now define several $p$-subgroups
$A$, $B_k$, $C$, $D$, and $E$ of $Aut(F_n)$.  Our goal is to
show that these are a complete listing of 
the distinct conjugacy classes of
$p$-subgroups of $Aut(F_n)$.

\begin{itemize}
\item There is an action of $\Z/p$ on the rose $R_{2(p-1)}$ given by
rotating the first $p$ leaves of the rose.  By looking at the
stabilizer of a marked graph with underlying graph $R_{2(p-1)}$, this
action gives us a subgroup $A \cong \Z/p$ of $Aut(F_{2(p-1)})$.  This
subgroup is a maximal $p$-subgroup, in the sense that no other
$p$-subgroup properly contains it.
\item For each $k \in \{0, \ldots, p-1 \}$, there is an action of
$\Z/p$ on $R_k \vee \Theta_{p-1} \vee R_{p-1-k}$ given by rotating the
edges in $\Theta_{p-1}$.  This action gives us a subgroup
$B_k \cong \Z/p$ of $Aut(F_{2(p-1)})$.  If $k \in \{1, \ldots, p-2
\}$, then $B_k$ is a maximal $p$-subgroup.
\item There is an action of $\Z/p$ on $\Phi_{2(p-1)}$ which gives us a
(non-maximal) $p$-subgroup $C \cong \Z/p$ of $Aut(F_{2(p-1)})$.
\item There is an action of $\Z/p \times \Z/p$ on $\Omega_{2(p-1)}$
given by having the first $\Z/p$ rotate one of the $\Theta$-graphs in
$\Omega_{2(p-1)}$ and having the second $\Z/p$ rotate the other
$\Theta$-graph.  This action gives us a subgroup $D \cong \Z/p \times
\Z/p$ of $Aut(F_{2(p-1)})$.  This subgroup is maximal among
$p$-subgroups, and it contains $B_0$, $B_{p-1}$, and $C$.
\item There is an action of $\Z/p \times \Z/p$ on $\Psi_{2(p-1)}$
given by having the first $\Z/p$ rotate one of the $\Theta$-graphs in
$\Psi_{2(p-1)}$ and having the second $\Z/p$ rotate the other
$\Theta$-graph.  This action gives us a subgroup $E \cong \Z/p \times
\Z/p$ of $Aut(F_{2(p-1)})$.  This subgroup is maximal among
$p$-subgroups, and it contains $B_{p-1}$ and $C$.
\end{itemize}

\begin{prop} \label{t2}
Every $p$-subgroup of $Aut(F_{2(p-1)})$
is conjugate to one of
$$A, B_0, \ldots , B_{p-1}, C, D, \hbox{ or } E.$$
\end{prop}

\begin{pf*}{Proof.} As we asserted earlier, every nontrivial 
$p$-subgroup $P$ is either $\Z/p$ or $\Z/p \times \Z/p$.
From Zimmerman's realization theorem, 
this subgroup is realized by an action on
a reduced basepointed marked graph $(\eta,\Gamma,*)$. 

If $P=\Z/p$, then Proposition \ref{t1} gives us that $\Gamma$ is 
one of
$R_n$, $R_k \vee \Theta_{p-1} \vee R_{p-1-k}$,
$\Psi_n$, or $\Omega_n$.  If $\Gamma$ is $R_n$ then $P$ is conjugate
to $A$.  If $\Gamma$ is $R_k \vee \Theta_{p-1} \vee R_{p-1-k}$ then
$P$ is conjugate to $B_k$.  Finally, note that
by collapsing different invariant forests the action
of $\Z/p$ on $\Phi_n$ gives a diagonal action of $\Z/p$ on
both $\Psi_n$ and $\Omega_n$.  Hence if $\Gamma$ is either
$\Psi_n$ or $\Omega_n$ then $P$ is conjugate to $C$.

Next, suppose $P=\Z/p \times \Z/p = (\alpha) \times (\beta)$.
The first cyclic summand must rotate $p$ edges
$e_1, e_2, \ldots, e_p$ of $\Gamma$. 
Without loss of
generality, we may assume that the basepoint $*$ is one of the
endpoints of each $e_i$.  Now if $\beta$ sends all of the $e_i$ to
another whole collection 
$\beta e_i$ (with $\{ e_i \}$ disjoint from $\{ \beta e_j \}$) then
the basepoint $*$ must be one of the endpoints of each $\beta^j e_i$ 
also; 
therefore, we obtain at least $p^2$ edges emanating from the 
basepoint $*$
which are moved by $\alpha$ and $\beta$.  This implies that the rank 
of
$\pi_1(\Gamma)$ is at least $p(p-1)$ (i.e., the best that can happen 
is
that $p$ copies of $\Theta_{p-1}$ are wedged together at the 
basepoint),
which is too large as $p \geq 3$.

So $\beta$ does not send the $e_i$ to another whole collection
$\beta e_i$ of edges disjoint from the $e_i$.
Without loss of generality, we may assume
$(\beta)$ fixes the edges $e_i$ (by replacing $\beta$ with
$\beta - \alpha^j$ if necessary.)  Hence the collection $\{ e_i \}$
is $P$-invariant. Now $\beta$
must rotate $p$ other edges $f_1, f_2, \ldots, f_p$.  
As $\Gamma$ is reduced, the $e_i$ do not form a subforest.
In addition, they do not form a $p$-gon as $*$ is an endpoint of 
each of 
them.  Hence the $e_i$ form either a rose or a $\Theta$-graph by 
the logic of 
Lemma \ref{tr1}.  They do not form a rose, else the existence of the 
edges $f_i$ forces the
rank of $\pi_1(\Gamma)$ to be at least $2p-1$.  Hence the $e_i$ 
form a 
$\Theta$-graph $\Theta_{p-1}$.  A similar argument shows that the 
$f_i$ 
also form a $\Theta_{p-1}$.
Hence $\Gamma$ is either $\Psi_n$ or
$\Omega_n$.  In the former case $P$ is conjugate to $E$ while
in the latter case it is conjugate to $D$. \END \end{pf*}

In order to determine whether or not two different graphs
give us the same elementary abelian $p$-subgroup,
we need to use work of Krstic in \cite{[K]} involving
{\em Nielsen transformations.}  More will be said on this
in the next section, but for now we briefly 
recall Krstic's definition of
Nielsen transformations and state the result of
Krstic that we need.

\begin{defn}[Nielsen transformation] \label{tr2}
Let $G$ be a finite subgroup of \newline
$Aut(F_n)$, which is realized by an
action of $G$ on a reduced, basepointed graph $\Gamma$.  
Let $V$ and $E$ be the vertex and oriented edge sets,
respectively, of $\Gamma$.  Finally, let $\iota, \tau : E \to V$ be
maps which give the initial and terminal points, respectively, of
oriented edges.
If there are two edges $e$ and $f$ of
$\Gamma$ such that:
\begin{itemize}
\item $e$ and $f$ are in different orbits (i.e., $f \not
\in e G \cup \bar e G$),
\item $\tau e = \tau f$, and
\item $stab(e) \subseteq stab(f)$.
\end{itemize}
then there is an \emph{admissible Nielsen transformation}
$<e,f>$ from $\Gamma$ to
a new graph $<e,f>\Gamma$.
The graph $\Gamma' = <e,f>\Gamma$ has the same vertex and
edge sets $V$ and $E$ as $\Gamma$; however, the map $\tau' : E \to V$
which gives the terminal point of an edge is changed as follows.  For
edges $h$ not in the orbit of $e$, set $\tau ' h=\tau h$; but set
$\tau' eg = \iota f g$ for $g \in G$.
\end{defn}

If a sequence of Nielsen transformations can change
a $G$-graph $\Gamma_1$ into a $G$-graph $\Gamma_2$,
then $\Gamma_1$ is said to be {\em Nielsen equivalent} to $\Gamma_2$.

We need the following theorem of Krstic, Theorem 2 from
\cite{[K]}:

\begin{thm}[Krstic] \label{tr4}
Let the graphs $\Gamma_1$ and $\Gamma_2$ realize the same subgroup
$G$ of
$Aut(F_n)$.  If $\Gamma_1$ and $\Gamma_2$ are reduced as
$G$-graphs then they are Nielsen equivalent,
up to an equivariant isomorphism (a basepoint preserving
isomorphism.)
\end{thm}

\begin{prop} \label{tr5}
The subgroups
$A$, $B_0, \ldots , B_{p-1}$, $C$, $D$, $E$
are in distinct conjugacy classes.  The diagram of
subgroups up to conjugacy is
\begin{equation} \label{e4}
\matrix{ & & B_0 & & \cr
           & & \downarrow & & \cr
           & & D & & \cr
           & \nearrow & & \nwarrow & \cr
           B_{p-1} & & & & C \cr
           & \searrow & & \swarrow & \cr
           & & E & & \cr}
\end{equation}
\end{prop}

\begin{pf*}{Proof.} We apply Theorem \ref{tr4}.
Any two reduced graphs realized by the same subgroup
of $Aut(F_n)$ can be connected up by a sequence of Nielsen
transformations.
This could only occur for two distinct graphs listed in
Proposition \ref{t1} when $P=\Z/p$ and the two graphs are
$\Psi_n$ and $\Omega_n$.  In
this case, $P$ is conjugate to the subgroup $
C$ of $Aut(F_n)$. \END \end{pf*}

\section{Normalizers of $p$-subgroups of $Aut(F_n)$} \label{c3}

The structure of centralizers of finite subgroups of $Aut(F_n)$ is
given by Krstic in \cite{[K]} in some detail.
He shows that, for a finite subgroup $G$ of $Aut(F_n)$, an element
of the centralizer $C_{Aut(F_n)}(G)$ is a product of
Nielsen transformations followed by a centralizing
graph isomorphism. Moreover, one of the
main propositions of his paper can be used to yield information about
the structure of normalizers of finite subgroups, which is what
we are interested in here.

We now recall some definitions and theorems from Krstic \cite{[K]}.

\begin{defn}[Nielsen isomorphism] \label{tr3}
Let $<e,f>$ be a Nielsen transformation from $\Gamma$ to $\Gamma'$
(see Definition \ref{tr2}.)
A \emph{Nielsen
isomorphism} is the isomorphism between fundamental groupoids:
$$<e,f> : \Pi(\Gamma) \to \Pi(\Gamma')$$
given by $<e,f>h=h$ if h is not in the orbit of $e$ and
$<e,f>ex=(ef)x$ for $x \in G$.
\end{defn}

The following is Proposition $4'$ from \cite{[K]}:

\begin{prop}[Krstic] \label{tr6}
If $\Gamma_1$ and
$\Gamma_2$ are finite basepointed, reduced 
 $G$- \newline 
graphs and
$$F : \Pi(\Gamma_1) \to \Pi(\Gamma_2)$$
is an equivariant groupoid isomorphism
which preserves the base vertex, then there
exists a product $T$ of Nielsen transformations and an equivariant
isomorphism of basepointed graphs
$$H: T\Gamma_1 \to \Gamma_2$$ such that
$F=HT$.
\end{prop}

Recall that $\pi_1(\Gamma_i)$ is a sub-groupoid of
$\Pi(\Gamma_i)$. Krstic uses the above theorem to
get information about maps between fundamental groups.
He shows the following, which is Corollary $1'$ in \cite{[K]}:

\begin{cor}[Krstic] \label{tr7}
Let $\Gamma_1$ and $\Gamma_2$ be reduced pointed $G$-graphs of
rank $\geq 2$.  Then every equivariant isomorphism
$$\pi_1(\Gamma_1,*) \to \pi_1(\Gamma_2,*)$$
is the restriction of an equivariant isomorphism
$$\Pi(\Gamma_1) \to \Pi(\Gamma_2).$$
\end{cor}

For a $G$-graph $\Gamma$ and an automorphism $\phi : G \to G$, let
$\phi(\Gamma)$ be the $G$-graph with underlying graph $\Gamma$ and
$G$-action given by $xg := x\phi(g)$ where $g \in P$, $x \in V(\Gamma)
\cup E(\Gamma)$, and the latter multiplication is given by the
standard action on $\Gamma$.  In other words, $\phi(\Gamma)$
is just the graph $\Gamma$ with the $G$-action ``twisted''
by $\phi$.

Realize the finite subgroup $G$ of $Aut(F_n)$ by a
marked graph $\eta : R_n \to \Gamma$ where
the induced action of $G$ on $\Gamma$
is reduced.

\begin{prop} \label{t3}
For every element $\alpha$ of
$N_{Aut(F_n)}(G)$ there exists an automorphism $\phi$ of
$G$ such that $\alpha$ is realized by
a $G$-equivariant isomorphism
$$\pi_1(\Gamma,*) \to \pi_1(\phi(\Gamma),*).$$
\end{prop}

\begin{pf*}{Sketch of proof.}
Define $N_{Aut(F_n)}'(G)$ to be the set of all
equivalence classes of pairs
$(\phi,\psi)$
where $\phi \in Aut(G)$ and
$\psi : \Gamma \to \Gamma$ is a basepoint preserving, continuous
surjection of graphs such that
$$\psi_\# : \pi_1(\Gamma,*) \to \pi_1(\phi(\Gamma),*)$$
is a $G$-equivariant group isomorphism.
Two such pairs $(\phi_1,\psi_1)$ and $(\phi_2,\psi_2)$
are {\em equivalent} if $\phi_1=\phi_2$ and
$$(\psi_1)_\# = (\psi_2)_\# : \pi_1(\Gamma,*) \to 
\pi_1(\phi(\Gamma),*)$$
on the level of fundamental groups.

Define a group operation {\em composition} in
$N_{Aut(F_n)}'(G)$ in the obvious way by letting
$$(\phi_2,\psi_2) \circ (\phi_1,\psi_1) =
(\phi_2 \circ \phi_1, \psi_2 \circ \psi_1).$$

To prove the proposition, it suffices to show that
there is a group isomorphism
$$\xi : N_{Aut(F_n)}(G) \to N_{Aut(F_n)}'(G).$$

This isomorphism is defined as follows.
Choose a
fixed homotopy inverse $\bar \eta : \Gamma \to R_n$; i.e., a map such
that $$(\eta \bar \eta)_\# = 1 \in Aut(\pi_1(\Gamma))$$ and
$$(\bar \eta \eta)_\# = 1 \in Aut(\pi_1(R_n)).$$
If $\alpha \in N_{Aut(F_n)}(G)$, then $\alpha$ induces an
automorphism $\phi$ of $G$ via conjugation: $\phi(g)=\alpha g
\alpha^{-1}$.  In addition, since $\alpha \in Aut(F_n)$ it
corresponds to a map $\tilde \alpha : R_n \to R_n$, allowing us to
set $\psi = \eta \tilde \alpha \bar \eta$. Now define $\xi$ by sending
$\alpha$ to the pair $(\phi,\psi)$.

We leave it as an exercise for the reader to verify that $\xi$ as
defined actually is a group isomorphism. \END \end{pf*}

As a result of Proposition \ref{t3}, we obtain the following:

\begin{prop} \label{t4} An element of $N_{Aut(F_n)}(G)$ is
a product of Nielsen transformations followed by a normalizing
graph isomorphism; that is, if $(\phi, \psi) \in N_{Aut(F_n)}'(G)$,
then there exist a product $T$ of Nielsen isomorphisms and 
a graph isomorphism $H$ such that $\psi_{\#} = HT$:
$$\psi_\#: \pi_1(\Gamma,*) \buildrel T \over \to
\pi_1(T\Gamma,*) \buildrel H \over \to
\pi_1(\phi(\Gamma),*).$$
\end{prop}

\begin{pf*}{Proof.} Any element $(\phi, \psi)$ of
$N_{Aut(F_n)}'(G)$ induces a $G$-equivariant map
$$\pi_1(\Gamma,*) \buildrel \psi_\# \over \to
\pi_1(\phi(\Gamma),*).$$
From Corollary \ref{tr7}, $\psi_\#$ is the restriction of
a $G$-equivariant map $F$ between fundamental groupoids:
$$\Pi(\Gamma) \buildrel F \over \to \Pi(\phi(\Gamma)).$$
Now from Proposition \ref{tr6} we know that there
exists a product $T$ of Nielsen transformations
starting with the graph $\Gamma$, and a
basepoint preserving graph isomorphism
$H : T\Gamma \to \phi(\Gamma)$ such that $F$ is the map
induced by $HT$: 
$$F: \Pi(\Gamma) \buildrel T \over \to
\Pi(T\Gamma) \buildrel H \over \to \Pi(\phi(\Gamma)).$$
By restricting the above map of fundamental groupoids to
one of fundamental groups, we have that $\psi_\#$ is
also the map induced by $HT$.
\END \end{pf*}

\bigskip

We now use Proposition \ref{t4} to calculate the normalizers
of the subgroups
$A$, $B_0, \ldots , B_{p-1}$, $D$, and $E$
listed in Proposition \ref{tr5}.  We do this so that their
cohomology groups will be easier to calculate in a later section.
We will not need
to calculate the normalizer of $C$ in this explicit manner,
and later will find its cohomology through more geometric means.

\begin{lemma} \label{tr8}
$$N_{Aut(F_n)}(A) \cong N_{\Sigma_p}(\Z/p) \times \bigl( (F_{p-2} 
\times F_{p-2})
\rtimes (\Z/2 \times Aut(F_{p-2}))  \bigr)$$
where $\Z/2$ acts by exchanging the two copies of $F_{p-2}$ and
$Aut(F_{p-2})$ acts diagonally on the two copies of $F_{p-2}$.
\end{lemma}

\begin{pf*}{Proof.} Write $F_n = \langle x_1, \ldots, x_n \rangle$ 
as the free 
group on the letters $x_j$.  In the above decomposition of 
$N_{Aut(F_n)}(A)$, the group  $N_{\Sigma_p}(\Z/p)$ corresponds to 
automorphisms of $F_n$ which permute the first $p$ letters
$x_1, \ldots, x_p$. In the decomposition, the $i$th generator of the
first copy of $F_{p-2}$ corresponds to an automorphism of $F_n$
which, for $j \leq p$, sends $x_j \mapsto x_j x_{p+i}^{-1}$ and
which fixes $x_j$ for $j > p$.  The $i$th generator of the
second copy of $F_{p-2}$ corresponds to an automorphism of $F_n$
which, for $j \leq p$, sends $x_j \mapsto x_{p+i} x_j$ and
which fixes $x_j$ for $j > p$.  Next, the $\Z/2$ comes from 
the map which, for $j \leq p$, sends $x_j \mapsto x_j^{-1}$
and which fixes the remaining $x_j$.  Finally, the $Aut(F_{p-2})$
comes from automorphisms which fix the first $p$ generators
$x_1, \ldots, x_p$ and which act on the latter generators
by identifying $F_{p-2}$ with 
$\langle x_{p+1}, x_{p+2}, \ldots, x_n \rangle.$ 

Recall that $A$ comes from the action of $\Z/p$ on the first $p$
petals of the rose $R_{2(p-1)}$.  From Proposition \ref{t4}
any element of $N_{Aut(F_n)}(A)$ is induced by a
product of Nielsen transformations $T$ followed by a
graph isomorphism
$$H: TR_{2(p-1)} \to \phi(R_{2(p-1)})$$
for some $\phi \in Aut(A)$.  The only types of Nielsen transformations
that are possible are
\begin{itemize}
\item Nielsen transformations obtained by pulling
either the front ends or the back ends of the first $p$ petals
of the rose uniformly around paths in the last $p-2$ petals.  This
subgroup is isomorphic to $F_{p-2} \times F_{p-2}$. 
\item Nielsen transformations contained
entirely in the copy of $Aut(F_{p-2})$ corresponding to
graph automorphisms and Nielsen transformations involving the
last $p-2$ petals.
\end{itemize}

Note that for any product $T$
of the above Nielsen transformations, the $A$-graph $TR_{2(p-1)}$ 
is exactly
the same as the $A$-graph $R_{2(p-1)}$. 
Hence from Proposition \ref{t4}, 
any element of $N_{Aut(F_n)}(A)$ is induced by a
product of Nielsen transformations $T$ followed by a
graph isomorphism
$$H: TR_{2(p-1)} = R_{2(p-1)} \to \phi(R_{2(p-1)})$$
for some $\phi \in Aut(A)$.
That is, the only ``normalizing graph automorphisms''
we need to examine are
automorphisms of one particular graph $R_{2(p-1)}$ which are
in the normalizer $N_{Aut(F_n)}(A)$.

The normalizing graph automorphisms are:
\begin{itemize}
\item Those involving just the last $p-2$ petals of the rose.
As the action of $A$ on these petals is trivial,
the normalizer of $A$ contains any graph automorphism
involving just those last $p-2$ petals.  All of these graph 
automorphisms
are contained in 
the copy of $Aut(F_{p-2})$ obtained from graph automorphisms and
Nielsen transformations involving the last $p-2$ petals.
\item Normalizing graph automorphisms of the first $p$
petals of the rose.  This gives a subgroup of $N_{Aut(F_n)}(A)$ 
that is
isomorphic to $N_{\Sigma_p}(\Z/p)$.
\end{itemize}

We leave it as an exercise for the reader to show that all of
the various subgroups of $N_{Aut(F_n)}(A)$ now fit together as 
described in the statement of the lemma. \END \end{pf*}

A final remark about the structure
of the subgroup $N_{Aut(F_n)}(A)$, described
above in the proof of Lemma \ref{tr8},
is appropriate here:
\begin{remark} \label{tr9} 
\end{remark}
Consider a subgroup $\langle \omega \rangle \cong \Z/2$ of 
$Aut(F_{p})$
corresponding to the action of $\Z/2$ on $R_p$
given by switching 
just the first two petals of the
rose.  Note that
$$N_{Aut(F_{p})}(\omega) = C_{Aut(F_{p})}(\omega)=
\Z/2 \times \bigl( (F_{p-2} \times F_{p-2})
\rtimes (\Z/2 \times Aut(F_{p-2}))  \bigr)$$
where the action of $\Z/2 \times Aut(F_{p-2})$
on $F_{p-2} \times F_{p-2}$ in the semidirect
product is the same as before (that is, as in $N_{Aut(F_n)}(A)$.)
Consequently, we see that
$$N_{\Sigma_p}(\Z/p) \times N_{Aut(F_p)}(\omega) \cong \Z/2 \times 
N_{Aut(F_n)}(A)$$
and hence
$$\hat H^*(N_{\Sigma_p}(\Z/p) \times N_{Aut(F_p)}(\omega); 
\Z_{(p)})
= \hat H^*(N_{Aut(F_n)}(A); \Z_{(p)})$$
because $p \geq 3$ and so the first summand
$\Z/2$ in $\Z/2 \times N_{Aut(F_n)}(A)$ can be ignored
when taking Farrell cohomology with $\Z_{(p)}$
coefficients.

\begin{lemma} \label{tr10}
$$N_{Aut(F_n)}(B_k) \cong N_{\Sigma_p}(\Z/p) \times (F_k \rtimes 
Aut(F_k))
\times Aut(F_{p-1-k})$$
\end{lemma}

\begin{pf*}{Proof.} Write $F_n = \langle x_1, \ldots, x_n \rangle$ 
as the free 
group on the letters $x_j$.  In the above decomposition of 
$N_{Aut(F_n)}(A)$, the group  $N_{\Sigma_p}(\Z/p)$ is obtained as
follows.  The permutation $(1 2 \ldots p)$ of $\Sigma_p$
corresponds to the automorphism of $F_n$ which fixes $x_j$ if
$j \leq k$, sends $x_j \mapsto x_{j+1}$ if
$k+1 \leq j \leq p+k-2$, sends 
$x_{p+k-1} \mapsto x_{p+k-1}^{-1} x_{p+k-2}^{-1} \ldots x_{k+1}^{-1}$,
and sends $x_j \mapsto x_{k+1}^{-1} x_j x_{k+1}$ if $j \geq p+k$.
The transposition $(1 2)$ of $\Sigma_p$ corresponds to the 
automorphism
which sends $x_{k+1} \mapsto x_{k+1}^{-1}$, 
sends $x_{k+2} \mapsto x_{k+1} x_{k+2}$,
sends $x_j \mapsto x_{k+1}^{-1} x_j x_{k+1}$ if $j \geq p+k$,
and fixes all other letters $x_j$.  Since $(1 2)$ and $(1 2 \ldots p)$
generate $\Sigma_p$, this suffices to define a copy of $\Sigma_p$ in
$Aut(F_n)$.  Now let $N_{\Sigma_p}(\Z/p)$ in the above
decomposition correspond to the normalizer
of $\Z/p \cong \langle (1 2 \ldots p) \rangle$ in this copy of the
symmetric group. 

The $F_k$ which is being acted upon in the semidirect product 
in the above decomposition has its $i$th generator
corresponding to the automorphism of $F_n$ which sends $x_j$ 
to its conjugate $x_i^{-1} x_j x_i$ if $j \geq k+1$ and fixes 
$x_j$ otherwise.  The $Aut(F_k)$ is the above decomposition
cooresponds to automorphisms of $F_n$ which fix the last
$n-k$ generators $x_{k+1}, \ldots, x_n$ and which act on the
first $k$ generators $x_1, \ldots, x_k$ as $Aut(F_k)$ indicates.

Finally, the $Aut(F_{p-1-k})$ in the above decomposition corresponds
to automorphism which fix the first $p-1+k$ generators of $F_n$ and
act on the last $p-1-k$ generators by identifying
$Aut(F_{p-1-k})$ with 
$Aut(\langle x_{p+k}, x_{p+k+1}, \ldots, x_n \rangle)$.

The remaining part of the proof is similar to that for 
Lemma \ref{tr8}, and we only sketch it here, leaving the details
for the reader to verify. Also, just for this proof, let $\Gamma$ 
denote the
$B_k$-graph $R_k \vee \Theta_{p-1} \vee R_{p-1-k}$.
Recall that $B_k$ comes from the action of $\Z/p$ on the 
$\Theta$-graph in
$\Gamma$.

The normalizer
$N_{Aut(F_n)}(B_k)$ contains four types of Nielsen
transformations: ones contained
in the subgroup $Aut(F_k)$ of the normalizer obtained from
taking graph automorphisms and Nielsen transformations of
the rose $R_k$,
ones contained
in the subgroup $Aut(F_{p-1-k})$ of the normalizer obtained from
taking graph automorphisms and Nielsen transformations of
the rose $R_{p-1-k}$,
redundant ones -- which can be ignored since they are already
included in the Nielsen transformations in $Aut(F_{p-1-k})$ -- 
obtained by pulling the back edges of the
$\Theta$-graph uniformly around the $p-1-k$ petals of the rose
$R_{p-1-k}$ on the right, and
ones corresponding to the $F_k$ in the decomposition in the statement
of the lemma which are obtained by pulling the front edges of the
$\Theta$-graph uniformly around the $k$ petals of the rose
$R_k$ on the left.

As was the case for Lemma \ref{tr8} above,
for any product $T$ of any of the above
types of Nielsen transformations,
the $B_k$-graph $T\Gamma$ is exactly the
same as the $B_k$-graph $\Gamma$.
The normalizing graph automorphisms take one of three forms:
automorphisms contained in the
subgroup $Aut(F_k)$ 
obtained from graph automorphisms and Nielsen transformations
involving just the rose $R_k$,
automorphisms contained in the
subgroup $Aut(F_{p-1-k})$, and graph 
automorphisms of the $\Theta_{p-1}$
sitting inside of $\Gamma$ which yield the
$N_{\Sigma_p}(\Z/p)$ in the decomposition of $N_{Aut(F_n)}(B_k)$. 
\END \end{pf*}

\begin{lemma} \label{tr11}
$$N_{Aut(F_n)}(D) = N_{\Sigma_p}(\Z/p) \times N_{\Sigma_p}(\Z/p).$$
$$N_{Aut(F_n)}(E) =
(N_{\Sigma_p}(\Z/p) \times N_{\Sigma_p}(\Z/p)) \rtimes \Z/2.$$
\end{lemma}

\begin{pf*}{Proof.} The groups
$N_{Aut(F_n)}(D)$ and $N_{Aut(F_n)}(E)$ are the easiest to calculate 
of
all of the normalizers.
There are no admissible Nielsen transformations of the
reduced $(\Z/p \times \Z/p)$-graphs $\Omega_n$ or $\Psi_n$.
It follows that $N_{Aut(F_n)}(D)$ and $N_{Aut(F_n)}(E)$ are just
finite groups consisting of
normalizing graph automorphisms of $\Omega_n$ and
$\Psi_n$, respectively.  So we need only examine the
graphs $\Omega_n$ or $\Psi_n$ and see which graph
automorphisms are in the normalizers of the
respective $(\Z/p \times \Z/p)$-actions.  Direct
examination of these graphs yields that the normalizers
are as claimed above. \END \end{pf*}

\section{Fixed point sets of the normalizers} \label{c8}

Let $G$ be a finite subgroup of $Aut(F_n)$, realized by a reduced
graph $\Gamma$ as in the previous section.  Define the
{\em fixed point subcomplex} $X_n^G$
of $G$ in the
spine $X_n$ of auter space by
$$X_n^G = \{x \in X_n : xg = x \hbox{ for all } g \in G\}.$$

From \cite{[K-V]} and Theorem 12.1 in Part III of \cite{[JD]}, we have

\begin{fact} \label{tr12}
The space $X_n^G$ is a contractible, finite-dimensional
complex and \newline
$N_{Aut(F_n)}(G)$ acts on it with finite quotient and
finite stabilizers.
\end{fact}

A few spaces related to fixed point subcomplexes will come up so
frequently in our work that we will give them special names,
$X_p^\omega$, $Q_p^\omega$, $\tilde X_m$, and $\tilde Q_m$.
We now define these spaces.

Recall from Krstic and Vogtmann \cite{[K-V]} that an edge of a reduced
$G$-graph $\Gamma$ is {\em inessential} if it is contained in
every maximal $G$-invariant forest in $\Gamma$,
and that there is an $N_{Aut(F_p)}(G)$-equivariant deformation
retract of $X_n^G$ obtained by collapsing all inessential
edges of marked graphs in $X_n^G$.

\begin{defn}[$X_p^\omega$ and $Q_p^\omega$] \label{tr17}
Let $\omega$ be the automorphism of $F_p$ defined by
interchanging the first two basis elements of $F_p$
(cf. Remark \ref{tr9}),
and let $X_p^{\langle \omega \rangle}$ be the subcomplex of
$X_p$ fixed by the subgroup $\langle \omega \rangle \cong \Z/2$.
Define $X_p^\omega$ to be the associated
$N_{Aut(F_p)}(\omega)$-equivariant deformation retract
of $X_p^{\langle \omega \rangle}$ and let $Q_p^\omega$
be the quotient of $X_p^\omega$ by the action of
$N_{Aut(F_p)}(\omega)$.
\end{defn}

From Remark $12.16$ of \cite{[JD]} we know that
the dimension of $X_p^\omega$ is $2p-4$.
In the next section, we will use $X_p^\omega$ to study 
the cohomology of $N_{Aut(F_n)}(A)$.
The rest of this section, however, will be devoted to
defining the space $\tilde X_m$, finding
an alternative characterization of this space,
and explaining how the space relates to
to the normalizer
$N_{Aut(F_n)}(B_m)$.

\begin{defn}[$\tilde X_m$ and $\tilde Q_m$] \label{tr18}
Let $m$ be a positive integer. Let $\alpha$ be the automorphism
of $F_{m+2}$ defined by:
$$\left\{ \matrix{
x_i \mapsto x_i \hbox{ for } i \leq m, \hfill \cr
x_{m+1} \mapsto x_{m+2} \hfill \cr
x_{m+2} \mapsto x_{m+2}^{-1} x_{m+1}^{-1}. \hfill \cr
} \right.$$
The subgroup $\mathcal{Q}$ generated by $\alpha$ has order 3,
and is realized on the graph $R_m \vee \Theta_2$ by rotating the
edges of $\Theta_2$ cyclically.  Define $\tilde X_m$ to be the fixed 
point set
$X_{m+2}^\mathcal{Q}$, and let $\tilde Q_m$ be the quotient of
$\tilde X_m$ by $N_{Aut(F_{m+2})}(\mathcal{Q})$.
For $m=0$, we set
$\tilde X_0 = \tilde Q_0 = \{\hbox{a point}\}$
as a notational convenience.
\end{defn}

Let $\tilde \Gamma$ be the graph $R_m \vee \Theta_{2}$,
where the basepoint of
the resulting graph $\tilde \Gamma$ is the center of the rose
$R_m$.
Label the $3$ edges of the $\Theta$-graph in $\tilde \Gamma$
as $e_1, e_2, e_3$ and orient them so that they begin
at the basepoint of $\tilde \Gamma$.
Construct a specific marked graph
$$\tilde \eta: R_{m+2} = R_{m} \vee R_{2} \to R_{m} \vee 
\Theta_{p-1}$$
by sending the first rose $R_{m}$ in $R_{m} \vee R_{2}$
to $R_{m}$ in $R_{m} \vee \Theta_{2}$ via the identity map.
Then send the second rose $R_{2}$ in $R_{m} \vee R_{2}$
to $\Theta_{2}$ in $R_{m} \vee \Theta_{2}$ by
sending the $i$th petal to $e_i * \bar e_{i+1}$.
As noted above in Definition \ref{tr18},
the subgroup $\mathcal{Q}$ acts on this marked graph by rotating the
edges of $\Theta_2$ cyclically.

We will
now study the structure of
$\tilde X_m$ and $\tilde Q_m$
in detail.
Recall from Lemma \ref{tr10} that
$$N_{Aut(F_n)}(B_k) \cong N_{\Sigma_p}(\Z/p) \times (F_k \rtimes 
Aut(F_k))
\times Aut(F_{p-1-k}).$$

In a similar way, we could use
Proposition \ref{t4} just as Lemma \ref{tr10} did to
obtain that
$$N_{Aut(F_{m+2})}(\mathcal{Q}) \cong \Sigma_3
\times (F_m \rtimes Aut(F_m)).$$
The $\Sigma_3 \cong N_{\Sigma_3}(\Z/3)$ comes from normalizing
graph automorphisms of the $\Theta_{2}$, the
$F_{m}$ comes from Nielsen moves done by pulling
the back ends of all of
the edges $e_i$ uniformly around some loop in the rose
$R_{m}$, and the $Aut(F_{m})$ comes from Nielsen
transformations and graph automorphisms concerning just
the petals of the rose $R_{m}$.
  
Let us examine the Nielsen transformations corresponding
to the $F_{m}$ in a little more detail.
Say the petals of the rose are $r_1, \ldots, r_{m}$.
The group
$F_{m}$ is generated by the $m$ Nielsen moves
$a_k := <\bar e_1, r_k>$ where $<\bar e_1, r_k>$
fixes the rose and sends
$e_j$ to $\bar r_k * e_j$ for all $j \in \{1, 2, 3\}$.
Observe that the 
marked graph $\tilde \eta \cdot a_k$
sends $R_{m}$ identically to
$R_{m}$ as before; however, the $i$th petal of the
second rose is now sent to $r_k * e_i * \bar e_{i+1}
* \bar r_k$ (recalling that we want a right action on marked
graphs and so $a_k$ acts on $\pi_1(\Gamma)$ by $a_k^{-1}$.)   
Hence some word
$$w=a_{k_1}, \ldots, a_{k_s} \in
F_m \subset N_{Aut(F_{m+2})}(\mathcal{Q})$$
in the letters $a_1, \ldots, a_{m}$ acts on the
marked graph $\tilde \eta$ by sending $R_{m}$
in $R_{m} \vee R_{2}$ identically
to the rose $R_{m}$ in $R_{m} \vee \Theta_{2}$;
however, the second rose $R_{2}$
is mapped to
$R_{m} \vee \Theta_{2}$ by conjugating the old
way it was mapped by the path
$r_{k_s} * \cdots * r_{k_1}$.

\begin{prop} \label{t11} The fixed point
space $\tilde X_{m}$ can be characterized
as the realization of the poset of equivalence classes of
pairs $(\alpha,f)$, where
$\alpha : R_{m} \to \Gamma_{m}$ is
a basepointed marked graph whose underlying graph $\Gamma_{m}$
has a special (possibly valence 2)
vertex which is designated as $\circ$,
$\circ$ may equal the basepoint $*$ of $\Gamma_{m}$,
and $f : I \to \Gamma_{m}$ is a homotopy
class (rel endpoints) of maps from
$*$ to $\circ$ in $\Gamma_{m}$. 
\end{prop}

\begin{pf*}{Proof.} By definition,
$\tilde X_{m}$ consists
of simplices in the spine $X_{m+2}$ that are
fixed by $\mathcal{Q}$.  It is the subcomplex
generated by marked graphs (i.e., vertices of $X_{m+2}$) which
realize the finite subgroup $\mathcal{Q} \subset Aut(F_{m+2})$.
From Theorem \ref{tr4}, any
two vertices in $\tilde X_{m}$ corresponding
to reduced marked graphs, are connected to
each other by Nielsen moves.  In other words,
they are connected to each other by Nielsen moves that
can be represented by
elements of
$$N_{Aut(F_{m+2})}(\mathcal{Q}) \cong \Sigma_3
\times (F_m \rtimes Aut(F_m)),$$
since the normalizer contains all of the relevant Nielsen
transformations.
In particular, the Nielsen moves come from
$<a_1, \ldots, a_{m}> \cong F_{m}$ and the
Nielsen moves in
$Aut(F_{m})$ involving only the petals of the rose
$R_{m}$ in $R_{m} \vee \Theta_{2}$.

Hence the reduced marked graphs representing vertices
of $\tilde X_{m}$ are of the form
$$\psi = \alpha \vee \beta :
R_{m} \vee R_{2} \to R_{m} \vee \Theta_{2}$$
where $\alpha : R_{m} \to R_{m}$
corresponds to any reduced marked graph representing
a vertex of $X_{m}$ and
$\beta : R_{2} \to R_{m} \vee \Theta_{2}$
sends the  $i$th petal to 
$$r_{k_s} * \cdots * r_{k_1} * e_i * \bar e_{i+1}
* \bar r_{k_1} * \cdots * \bar r_{k_s},$$
and corresponds to some word
$a_1, \ldots, a_s$ in $F_{m}$.
We could thus represent the reduced
marked graph $\psi$ more compactly as a pair
$(\alpha,f)$ where $\alpha: R_{m} \to R_{m}$
is any reduced marked graph and $f:I=[0,1] \to R_{m}$
is a homotopy class (rel the endpoints $f(0)=f(1)=*$) of maps
representing a path $r_{k_s} * \cdots * r_{k_1}$
in the rose $R_{m}$.

By considering
the stars, in $\tilde X_{m}$, of reduced
marked graphs
$(\alpha,f): R_{m+2} \to R_{m} \vee \Theta_{2}$,
we will 
obtain a characterization of all
marked graphs representing vertices of
$\tilde X_{m}$.
Recall that $\mathcal{Q}$ rotates the edges of $\Theta_{2}$
and leaves everything else fixed.

\begin{defn} \label{tr22} Let $G$ be a finite
subgroup of $Aut(F_r)$ for some integer $r$.
A marked graph
$$\eta^1 : R_{r} \to \Gamma^1$$
is a {\em $G$-equivariant blowup in the fixed point space
$X_r^G$}
of a marked graph
$$\eta^2 : R_{r} \to \Gamma^2$$
if there is a $1$-simplex $\eta^1 > \eta^2$ in $X_r^G$.
In other words, this happens exactly when we can
collapse a $G$-invariant forest in $\Gamma^1$
to obtain $\Gamma^2$. 
We often abbreviate this and
just say that $\eta_1$ is a {\em blowup} of $\eta_2$ or
that $\eta_2$ can be {\em blown up} to get $\eta_1$.
If the forest that we collapsed was just a tree, we
say that we {\em blew up} the corresponding vertex
of $\eta_2$ to get $\eta_1$.  
\end{defn}

\begin{claim} \label{tr23}
Any equivariant blowup of $\Gamma$ is a blowup of $R_m$,
with $\Theta_2$ attached at a point.
\end{claim}

\begin{pf*}{Proof.}
Suppose $\Gamma$ is obtained from $\Gamma'$ by collapsing an 
invariant forest
$F$.  It suffices to show that $F$ is fixed by the action of 
$\mathcal{Q}$, since then the
initial (resp. terminal) vertices of the edges mapping to $e_i$ 
must be the same.

Suppose $F$ is not fixed by $\mathcal{Q}$, and let $F_0$
be the union of edges of $F$ with non-trivial orbits.  Let $v$ be 
a terminal vertex of $F_0$, and let $e$ and $f$ be edges of $\Gamma'$
terminating at $v$ which are not in $F_0$.  If $e$ or $f$ is fixed by 
$\mathcal{Q}$, then $v$ is fixed; in particular, if $e$ or $f$ is in
$F-F_0$, then $v$ is fixed by $\mathcal{Q}$.  If $e$ and $f$ are not 
fixed, then since they are not in $F$, they must map to edges $e_i$ 
and $e_j$ of 
$\Gamma$, with $i \not = j$.  Since the map is equivariant, some 
element of $\mathcal{Q}$
takes $e$ to $f$.  But this element must then fix their common 
vertex, $v$, so
in all cases $v$ must be fixed by $\mathcal{Q}$.  Since all 
terminal vertices of
$F_0$ are fixed, $F_0$ must be fixed, contradicting the 
definition of $F_0$. \END \end{pf*}

Hence the only way to blow up $(\alpha,f)$ is
by blowing up the $R_{m}$ part of
$R_{m} \vee \Theta_{2}$. 
Think of the resulting marked graph as
some $(\hat \alpha, \hat f)$ where
$$\hat \alpha : R_{m} \to \Gamma_{m}$$ is any
marked graph, except that
the underlying graph $\Gamma_{m}$ has
one extra distinguished vertex, which we will call $\circ$
and which might have valence $2$,
aside from the basepoint $*$, 
and where
$$\hat f : I \to \Gamma_{m}$$
is a homotopy class of
paths from $*$ to $\circ$ in $\Gamma_{m}$.
We allow the possibility that $\circ$ might just be $*$.
The pair $(\hat \alpha, \hat f)$ is really
representing a marked
graph $R_{m+2} \to \Gamma_{m} \vee \Theta_{2}$
as follows.  The wedge $\vee$ connects the
point $\circ$ of $\Gamma_{m}$ to the left
hand vertex of $\Theta_{2}$.  The basepoint of
the resulting graph
$\Gamma_{m} \vee \Theta_{2}$
is whatever the old basepoint $*$ of $\Gamma_{m}$
was.  The first $m$ petals of
$R_{m+2}$ map to $\Gamma_{m}$ via $\hat \alpha$
and the $(m+i)$-th petal of $R_{m+2}$
first goes around the image of $\hat f$ from
$*$ to $\circ$, then goes around $e_i * \bar e_{i+1}$
of the $\Theta$-graph, and then finally
goes back around the image of $\hat f$ in
reverse from $\circ$ to $*$.

For $(\hat \alpha, \hat f)$ to be in the
star of $(\alpha,f)$ it has to be
the case that when the graph
$\Gamma_{m} \vee \Theta_{2}$
is collapsed to $R_{m} \vee \Theta_{2}$
the marking $\hat \alpha$ collapses
to $\alpha$ and the homotopy class
(rel endpoints) of maps $\hat f$ collapses
to $f$. This concludes the proof of
Proposition \ref{t11}. \END \end{pf*}

As in \cite{[C-V]}, there is an obvious
definition of when two of the marked graphs described in
Proposition \ref{t11} 
are equivalent. 

\begin{defn} \label{tr24}
Marked graphs
$(\alpha^1,f^1)$ and $(\alpha^2,f^2)$
are {\em equivalent} 
if there
is a homeomorphism $h$ from 
$\Gamma_{m}^1$ to $\Gamma_{m}^2$
which sends $*$ to $*$, $\circ$ to $\circ$,
such that
$$(h \alpha^1)_\# = (\alpha^2)_\# : \pi_1(R_m,*) \to 
\pi_1(\Gamma_m^2)$$
and such that the paths
$$ hf_1, f_2 : I \to \Gamma_m^2$$
are homotopic rel endpoints.
\end{defn}

By recalling how
$N_{Aut(F_{m+2})}(\mathcal{Q})$ acts on reduced marked
graphs and looking at $(\alpha,f)$ as being a
marked graph
$$R_{m+2} \to \Gamma_{m} \vee \Theta_{2}$$
in the star of some reduced marked graph
$$R_{m+2} \to R_{m} \vee \Theta_{2},$$ 
it is
direct to prove 

\begin{prop} \label{t12} The group
$$N_{Aut(F_{m+2})}(\mathcal{Q}) \cong
\Sigma_3 \times (F_{m} \rtimes Aut(F_{m}))$$
acts on a marked graph
$$\matrix{\hfill (\alpha,f): R_{m} \coprod I \to \Gamma_{m}\hfill 
\cr}$$
in $\tilde X_{m}$ as follows.
The subgroup $\Sigma_3$ 
permutes the edges of the $\Theta$ graph attached at $\circ$,
giving a marked graph which is equivalent to the
original one.  An element $a_{k_1} a_{k_2} \ldots a_{k_s}$
in $F_{m} = <a_1, \ldots, a_{m}>$ doesn't change
$\alpha$ at all, but sends the path $f$ to the
path 
$$\alpha(r_{k_s}) * \alpha(r_{k_{s-1}}) * \cdots * \alpha(r_{1}) * f$$
where $r_i$ is the $i$th petal of the rose $R_{m}$
in the domain of $\alpha$.
Lastly, an element $\phi \in Aut(F_{m})$
does not change $f$ at all and acts on
$\alpha: R_{m} \to \Gamma_{m}$
by precomposition:
$$(\alpha,f) \cdot \phi = (\alpha \circ \phi, f).$$
\end{prop}

Just as Proposition \ref{t11} gives us
a simple characterization of $\tilde X_m$, the following remark
provides a nice characterization of $\tilde Q_m$:

\begin{remark} \label{tr25} The quotient
space $\tilde Q_{m}$ of $\tilde X_m$
by $N_{Aut(F_{m+2})}(\mathcal{Q})$
can be characterized
as the realization of the poset of equivalence classes of
basepointed graphs $\Gamma_{m}$
which have 
a special (possibly valence 2)
vertex designated as $\circ$,
which may equal the basepoint $*$.
Two such graphs $\Gamma_{m}^1$ and $\Gamma_{m}^2$
are {\em equivalent} if there 
is a homeomorphism
$$h: \Gamma_{m}^1 \to \Gamma_{m}^2$$
such that $h(*)=*$ and $h(\circ)=\circ$.
Define the poset structure
on these graphs by
forest collapses.  That is,
$\Gamma_{m}^1 > \Gamma_{m}^2$
if there is a simplicial map
$g: \Gamma_{m}^1 \to \Gamma_{m}^2$
such that $g(*)=*$, $g(\circ)=\circ$,
and $g^{-1}(\hbox{vertices of }\Gamma_{m}^2)$ is
a subforest of $\Gamma_{m}^1$.
\end{remark}

In other words, the
quotient $\tilde Q_{m}$ is just the ``moduli space
of unmarked graphs $\Gamma_{m}$ with
$\pi_1(\Gamma_{m}) \cong F_{m}$
and where $\Gamma_{m}$ has two distinguished points.''

Recall from \cite{[H-V]} that the spine $X_m$ of auter space is a
deformation retraction of auter space $\A_m$.  Similarly, we can
construct a space $\tilde \A_{m}$ which deformation retracts
to $\tilde X_m$.  We can then think of $\tilde X_m$ as
being the ``spine'' of $\tilde \A_m$.

\begin{defn}[$\tilde \A_{m}$] \label{tr26}
Construct an
analog $\tilde \A_{m}$ of auter space
for $N(\mathcal{Q})$ by considering
markings
$$\matrix{\hfill (\alpha,f): R_{m} \coprod I \to \Gamma_{m}\hfill 
\cr}$$
where
the edges of $\Gamma_{m}$ are assigned
lengths which must sum to 1.  Just as
in \cite{[C-V]}, the space $\tilde \A_{m}$
deformation retracts to its spine
$\tilde X_{m}$.
\end{defn}

Now from \cite{[H-V]} we have
$dim(X_{m}) = dim(Q_{m}) = 2m-2$
and $dim(\A_{m}) = 3m-3$.  As the graph
$\Gamma_{m}$ in a particular marked graph
has possibly one extra vertex
$\circ$ of valence $2$, we see that
$dim(\tilde X_{m}) = dim(\tilde Q_{m}) = 2m-1$
and $dim(\tilde \A_{m}) = 3m-2$.

Now that we are more familiar with the
structure of our spaces
$\tilde X_m$ and $\tilde Q_m$,
we can proceed with describing how they are
related to the normalizers
$N_{Aut(F_n)}(B_k)$.

\begin{defn}[Action of $N_{Aut(F_n)}(B_k)$ on $\tilde X_k$] 
\label{tr20}
For $k \in \{0, \ldots, p-1\}$,
the map from
$$N_{Aut(F_n)}(B_k) \cong N_{\Sigma_p}(\Z/p) \times (F_k 
\rtimes Aut(F_k))
\times Aut(F_{p-1-k})$$
to
$$N_{Aut(F_{m+2})}(\mathcal{Q}) \cong \Sigma_3
\times (F_m \rtimes Aut(F_m))$$
given by projection to the second factor followed by
inclusion induces an action of 
$N_{Aut(F_n)}(B_k)$ on $\tilde X_k$.
Since $N_{\Sigma_p}(\Z/p) \subset N_{Aut(F_n)}(B_k)$ acts
trivially, the quotient of this action is $\tilde Q_k$. 
\end{defn}

We can also define an action of $N_{Aut(F_n)}(B_k)$
on $X_{p-1-k}$:

\begin{defn}[Action of $N_{Aut(F_n)}(B_k)$ on $X_{p-1-k}$] 
\label{tr21}
For $k \in \{0, \ldots, p-1\}$,
the map from
$$N_{Aut(F_n)}(B_k) \cong N_{\Sigma_p}(\Z/p) \times (F_k 
\rtimes Aut(F_k))
\times Aut(F_{p-1-k})$$
to
$Aut(F_{p-1-k})$
given by projection on the third factor induces an action of
$N_{Aut(F_n)}(B_k)$ on $X_{p-1-k}$,
with quotient $Q_{p-1-k}$.
\end{defn}

From Definition \ref{tr20} and Definition \ref{tr21},
there is an induced action of \newline
$N_{Aut(F_n)}(B_k)$ on
$\tilde X_k \times X_{p-1-k}$.
Note that the space
$\tilde X_k \times X_{p-1-k}$
is a product of posets and can be given a poset
structure by saying that
$$(\alpha^1,f^1) \times \beta^1 \geq
(\alpha^2,f^2) \times \beta^2$$
if
$$(\alpha^1,f^1) \geq (\alpha^2,f^2)
\hbox{ and }
\beta^1 \geq \beta^2.$$  This gives
$\tilde X_k \times X_{p-1-k}$
a simplicial structure.  Often we will
not use this simplicial structure, however,
and just use the cellular structure that comes
from products of simplices in
$\tilde X_k$ and $X_{p-1-k}$.

Because $\tilde X_k$ and $X_{p-1-k}$ are contractible from Fact 
\ref{tr21}
and \cite{[H-V]}, respectively, we see that

\begin{thm} \label{tr28}
The above induced action of
$N_{Aut(F_n)}(B_k)$ on the contractible
space $\tilde X_k \times X_{p-1-k}$
has finite stabilizers and
quotient
$\tilde Q_k \times Q_{p-1-k}$.
\end{thm}

It is interesting to note that the space
$\tilde X_k \times X_{p-1-k}$ is in fact
homeomorphic to the fixed point subcomplex
$X_n^{B_k}$.  We leave this as a
straightforward exercise for the
reader, as this fact will not be used in this paper.

\section{The cohomology of the normalizers} \label{c4}

In this section, we will compute the cohomology of the normalizers
of the subgroups
$$A, B_k, C, D, \hbox{ and } E$$
listed in
Proposition \ref{tr5}.
First, we list some helpful facts that will allow us to
compute these cohomology groups.

The following theorem of Swan's
(see \cite{[S]} or \cite{[A-M]}) is a standard
tool for computing the cohomology of groups:
\begin{thm}[Swan] \label{tr30}
If $G$ is a finite group with a $p$-Sylow
subgroup $P$ that is abelian, then
$$H^*(G; \Z_{(p)}) = H^*(P; \Z_{(p)})^{N_{G}(P)}.$$
\end{thm}

In part II of \cite{[JD]} we showed several low dimensional 
cohomology groups
having to do with $Aut(F_n)$ are zero.  The methods used are variants 
of those used by Hatcher and Vogtmann in \cite{[H-V]} and 
\cite{[V]} where they use a ``degree theorem'' to reduce complicated
calculations involving $X_n$ to more manageable calculations involving
subcomplexes of $X_n$.  These results are used in this paper to prove
that certain spectral sequences converge at the $E_2$-page.
We list the results here in the following fact:

\begin{fact} \label{p1t1}
\begin{enumerate}
\item $H^1(Q_p^\omega; \Z/p) = 0.$
\item $H^1(Q_m; \Z/p) = 0$ for $1 \leq m \leq p-1$.
\item $H^2(Q_{p-1}; \Z/p) = 0.$
\item $H^1(\tilde Q_m; \Z/p) = 0$ for $1 \leq m \leq p-1$.
\item $H^2(\tilde Q_{p-1}; \Z/p) = 0.$
\end{enumerate}
\end{fact}

We can now calculate the cohomology of the
normalizers of the various subgroups $A$, $B_k$, $C$,
$D$, and $E$.

\begin{lemma} \label{t5} $$\hat H^t(N_{Aut(F_n)}(A); \Z_{(p)})
\cong \left\{\matrix{
\Z/p \hfill &t \equiv 0  \hbox{ } (\hbox{mod } n) \hfill \cr
0 \hfill &t \equiv \pm 1  \hbox{ } (\hbox{mod } n) \hfill \cr
H^r(Q_{p}^\omega; \Z/p) \hfill 
&t \equiv r  \hbox{ } (\hbox{mod } n), \hfill \cr
 &2 \leq r \leq n-2 \hfill \cr
} \right.$$
\end{lemma}

\begin{pf*}{Proof.} 
We now define an action of $N_{\Sigma_p}(\Z/p) \times 
N_{Aut(F_p)}(\omega)$
on the space $X_{p}^\omega$ by
stipulating that $N_{\Sigma_p}(\Z/p)$ acts trivially and
that $N_{Aut(F_p)}(\omega)$ acts in the usual manner on 
$X_{p}^\omega$.
This in turn defines an action of $N_{Aut(F_n)}(A)$
on the contractible space $X_p^\omega$.  This action has finite
stabilizers and quotient $Q_p^\omega$.

The equivariant cohomology spectral sequence
for this action is
\begin{equation} \label{e5}
\matrix{
E_1^{r,s} = \prod_{[\delta] \in \Delta_{p-1}^r}
\hat H^s(stab_{N_{Aut(F_n)}(A)}(\delta); \Z_{(p)}) 
\hspace*{4cm} \hfill \cr
\hfill
\Rightarrow
\hat H^{r+s}(N_{Aut(F_n)}(A)); \Z_{(p)}) \cr}
\end{equation}
where $[\delta]$ ranges over the set $\Delta_{p-1}^r$ of orbits
of $r$-simplices $\delta$ in $X_{p}^\omega$.

We claim
that each $stab_{N_{Aut(F_n)}(A)}(\delta)$ is the direct sum of
$N_{\Sigma_p}(\Z/p)$ with a finite subgroup of
$N_{Aut(F_p)}(\omega)$ that does
not have any $p$-torsion.  This is because
Glover and Mislin showed in \cite{[G-M]} that
the only $p$-torsion in
$Aut(F_{p})$ comes from stabilizers of marked graphs with underlying
graphs $\Theta_{p}$, $\Theta_{p-1} \vee R_1$,
$\Xi_p$, or $R_p$, where the graph $\Xi_p$ is defined to be the
$1$-skeleton of the
cone over a $p$-gon.  Since $\omega$ acts by switching the
first two petals of the rose $R_p$,
$X_{p}^\omega$ does not contain any marked
graphs with underlying graph
$\Theta_{p}$, $\Theta_{p-1} \vee R_1$, or $\Xi_p$.
Although $X_{p}^\omega$ obviously does contain marked
graphs with underlying graph $R_p$, these will not worry us
as $p$ does not divide the orders of the stabilizers
-- under the action of
just $N_{Aut(F_p)}(\omega)$ -- of such marked
graphs.  For example,
$N_{Aut(F_p)}(\omega)$
does not contain
the permutation $(1 2 \ldots p)$ that
rotates the petals of the rose $R_p$
because
$$(1 2 \ldots p) \circ \omega \circ (p  \ldots 2 1)$$
is the permutation
$(2 3)$, which is not equal to $(1)$ or
$\omega = (1 2)$.

Thus for every $[\delta]$, we have
$$\hat H^s(stab_{N_{Aut(F_n)}(A)}(\delta); \Z_{(p)}) =
\hat H^s(N_{\Sigma_p}(\Z/p); \Z_{(p)})
= \hat H^s(\Sigma_p; \Z_{(p)}).$$
The $E_1^{r,s}$-page of the spectral sequence is $0$ in the rows
where $s \not = kn$ and a copy
of the cellular cochain complex with $\Z/p$-coefficients of the
$(n-2)$-dimensional complex $Q_{p}^\omega$ in rows $kn$.  
It follows that
the $E_2$-page has the form:
$$E_2^{r,s} = \left\{\matrix{
\Z/p \hfill &r=0 \hbox{ and } s=kn \hfill \cr
H^r(Q_{p}^\omega; \Z/p) \hfill &1 \leq r \leq n-2
\hbox{ and } s=kn \hfill \cr
0 \hfill &\hbox{otherwise } \hfill \cr} \right.$$
Hence we see that the spectral sequence converges at the $E_2$-page.

That $H^1(Q_{p}^\omega; \Z/p) = 0$ follows from 
part 1 of Fact \ref{p1t1}. \END \end{pf*}

\begin{lemma} \label{t6} 
$$\matrix{
\hfill \cr
\hat H^t(N_{Aut(F_n)}(B_0); \Z_{(p)}) \cong
\left\{\matrix{
\Z/p^2 \hfill &t = 0 \hfill \cr
(j+1) \Z/p \hfill &|t| = nj \not = 0 \hfill \cr
(j-1) \Z/p \hfill &|t| = nj-1 \hfill \cr
0 \hfill &|t| \equiv 1  \hbox{ } (\hbox{mod } n) \hfill \cr
0 \hfill &t \equiv 2  \hbox{ } (\hbox{mod } n) \hfill \cr
H^r(Q_{p-1};\Z/p) \hfill
&t \equiv r  \hbox{ } (\hbox{mod } n), \hfill \cr
 \hfill &3 \leq r \leq n-2 \hfill \cr
}\right.
\hfill \cr
\hfill & \hfill \cr
}$$
$$\matrix{
\hbox{For $k \in \{1,\ldots,p-2\},$} \hfill \cr
\hfill \cr
\hat H^t(N_{Aut(F_n)}(B_k); \Z_{(p)}) \cong
\left\{\matrix{
\Z/p \hfill &t \equiv 0  \hbox{ } (\hbox{mod } n) \hfill \cr
0 \hfill &t \equiv \pm 1, -2,  \hbox{ } (\hbox{mod } n) \hfill \cr
H^r(\tilde Q_{k} \times Q_{p-1-k};\Z/p) \hfill
&t \equiv r  \hbox{ } (\hbox{mod } n), \hfill \cr
\hfill &2 \leq r \leq n-3 \hfill \cr
} \right.
\hfill \cr
\hfill & \hfill \cr}$$
$$\matrix{
\hfill \cr
\hat H^t(N_{Aut(F_n)}(B_{p-1}); \Z_{(p)}) \cong
\left\{\matrix{
\Z/p^2 \oplus \Z/p \hfill &t = 0 \hfill \cr
(2j+1) \Z/p \hfill &|t| = nj \not = 0 \hfill \cr
0 \hfill &t = nj + 1 > 0 \hfill \cr
(2j-2) \Z/p \hfill &t = -nj+1 < 0 \hfill \cr
(2j-2) \Z/p \oplus  \hfill
&t = nj-1 > 0 \hfill \cr
\hfill H^{n-1}(\tilde Q_{p-1}; \Z/p) & \hfill \cr
H^{n-1}(\tilde Q_{p-1}; \Z/p) \hfill &t = -nj-1 < 0 \hfill \cr
0 \hfill &t \equiv 2  \hbox{ } (\hbox{mod } n) \hfill \cr
H^r(\tilde Q_{p-1}; \Z/p) \hfill &t \equiv r  \hbox{ } 
(\hbox{mod } n), \hfill \cr
 &3 \leq r \leq n-2 \hfill \cr
} \right.
\hfill \cr
\hfill \cr}$$
\end{lemma}

\begin{pf*}{Proof.} We use the equivariant cohomology spectral
sequence (\ref{e2}) corresponding to the action,
defined in Definition \ref{tr20} and Definition \ref{tr21},
of
$$N_{Aut(F_n)}(B_k)
\cong N_{\Sigma_p}(\Z/p) \times (F_k \rtimes Aut(F_k))
\times Aut(F_{p-1-k})$$
on the space $\tilde X_k \times X_{p-1-k}$.
From Theorem \ref{tr28}
the action of
$N_{Aut(F_n)}(B_k)$ on the contractible
space $\tilde X_k \times X_{p-1-k}$
has finite stabilizers and
quotient
$\tilde Q_k \times Q_{p-1-k}$.

Applying the spectral sequence (\ref{e2}) we obtain
the following $E_1$-page:
\begin{equation} \label{e6}
\matrix{
E_1^{r,s} = \prod_{[\delta] \in \Delta^r}
\hat H^s(stab_{N_{Aut(F_n)}(B_k)}(\delta); \Z_{(p)}) 
\hspace*{4cm} \hfill \cr
\hfill
\Rightarrow
\hat H^{r+s}(N_{Aut(F_n)}(B_k)); \Z_{(p)}) \cr}
\end{equation}
where $[\delta]$ ranges over the set $\Delta^r$ of orbits
of $r$-simplices $\delta$ in $\tilde X_k \times X_{p-1-k}$.
Vertices in $\Delta^0$ are
pairs of unmarked graphs 
$(\Gamma_1,\Gamma_2)$ where $\pi_1(\Gamma_1) \cong F_k$,
$\pi_1(\Gamma_2) \cong F_{p-1-k}$, the graph $\Gamma_1$ has
two distinguished points $*$ and $\circ$,
and the graph $\Gamma_2$ has one distinguished point $*$.

The stabilizer of this vertex of $\Delta^0$
under the action of
$$N_{Aut(F_n)}(B_k)
\cong N_{\Sigma_p}(\Z/p) \times (F_k \rtimes Aut(F_k))
\times Aut(F_{p-1-k})$$
is isomorphic to
$$N_{\Sigma_p}(\Z/p) \times Aut(\Gamma_1) \times Aut(\Gamma_2)$$
where by $Aut(\Gamma_1)$ we mean graph automorphisms of
$\Gamma_1$ that preserve
both distinguished points, and $Aut(\Gamma_2)$ is the group of graph
automorphisms of $\Gamma_2$ that preserve its distinguished point.

There are three cases, depending upon what value
$k$ takes.

\smallskip
\noindent CASE 1:
$k \in \{1,\ldots,p-2\}.$
Consider an $r$-simplex $\delta$
$$((\alpha^0,f^0),\beta^0) > ((\alpha^1,f^1),\beta^1) >
\cdots > ((\alpha^r,f^r),\beta^r)$$
of $\tilde X_k \times X_{p-1-k}$.
Let $\Gamma^1$ be the underlying graph of
$(\alpha^0,f^0)$ and let $\Gamma^2$ be the underlying
graph of $\beta^0$.
The stabilizer of $\delta$ is a subgroup of the group
$$N_{\Sigma_p}(\Z/p) \times Aut(\Gamma_1) \times Aut(\Gamma_2).$$
The finite group $Aut(\Gamma_1)$
has no $p$-torsion
since $k < p-1$ and so none of the
underlying graphs of
marked graphs in $\tilde X_k$ have any $p$-symmetry.
Similarly,
the finite group $Aut(\Gamma_2)$
has no $p$-torsion since $k > 0$ and so none of the
underlying graphs of
marked graphs in $X_{p-1-k}$ have any $p$-symmetry.
Thus we have that
$$\hat H^*(stab(\delta); \Z_{(p)}) \cong
\hat H^*(N_{\Sigma_p}(\Z/p); \Z_{(p)}) \cong
\hat H^*(\Sigma_p; \Z_{(p)}).$$

Since the above holds for every simplex $\delta$,
we see that
the spectral sequence (\ref{e6}) has $E_2$ page

$$E_2^{r,s} = \left\{\matrix{
\Z/p \hfill &r=0 \hbox{ and } s=kn \hfill \cr
H^r(\tilde Q_{k} \times Q_{p-1-k}; \Z/p) \hfill &1 \leq r \leq n-3
\hbox{ and } s=kn \hfill \cr
0 \hfill &\hbox{otherwise } \hfill \cr} \right.$$
because the dimension of $\tilde Q_{k} \times Q_{p-1-k}$
is $2p-5$.

Now apply parts 2 and 4 of Fact \ref{p1t1} 
 to conclude 
that all of the groups
$H^1(\tilde Q_{k} \times Q_{p-1-k};\Z/p)$
are zero.
The lemma now follows for $k \in \{1,\ldots,p-2\}.$

\smallskip
\noindent CASE 2:
$k = 0.$  Then the simplices $\delta$
in spectral sequence (\ref{e6})
are all in $Q_{p-1}$.
Since only one graph in $Q_{p-1}$ has
$p$-symmetry, namely the graph $\Theta_{p-1}$,
we have

$$stab_{N_{Aut(F_n)}(B_k)}(\delta) = \left\{\matrix{
N_{\Sigma_p}(\Z/p) \times \Sigma_p \hfill
&\hbox{if } \delta \hbox{ has underlying graph } \hfill \cr
\hfill & \hfill \Theta_{p-1} \cr
\hfill & \hfill \cr
N_{\Sigma_p}(\Z/p) \times H \hfill &\hbox{otherwise, where 
$H$ is a} \hfill \cr
\hfill &\hfill \hbox{group with } p \not | \hbox{ } |H| \cr} \right.$$

Arguments similar to those in Lemma \ref{t5} show that the $E_2$ 
page has
the form

$$E_2^{r,s} = \left\{\matrix{
\Z/p^2 \oplus \Z/p \hfill & r=0 \hbox{ and } s=0 \hfill \cr
H^r(Q_{p-1}; \Z/p) \hfill &1 \leq r \leq n-2, s = kn \hfill \cr
(j-1)\Z/p \hfill &r=0 \hbox{ and } |s|=nj-1 \hfill \cr
(j+1)\Z/p \hfill &r=0 \hbox{ and } |s|=nj \not = 0  \hfill \cr
0 \hfill &\hbox{otherwise } \hfill \cr} \right.$$

Because part 3 of Fact \ref{p1t1}  gives that
$H^2(Q_{p-1}; \Z/p) = 0$,
the
differentials $d_2: E_2^{0,-nj+1} \to E_2^{2,-nj}$ are zero
and we see that the spectral sequence converges at the $E_2$ page.
Lastly, by part 2 of Fact \ref{p1t1}, $H^1(Q_{p-1}; \Z/p)=0$
and we are done with the case $k=0$.

\smallskip
\noindent CASE 3:
$k = p-1.$
Then the simplices $\delta$
in spectral sequence (\ref{e6})
are all in $\tilde Q_{p-1}$.
Now only two graphs in $\tilde Q_{p-1}$ have
$p$-symmetry.
One is the graph $\Theta_{p-1}^1$
where both $*$ and $\circ$
are the left hand vertex of the
$\Theta$-graph.  The other is
the graph $\Theta_{p-1}^2$
where
$*$ is the vertex on the left side of the $\Theta$-graph
and $\circ$ is the vertex on the right side.
Each of these graphs gives a vertex of $\tilde Q_{p-1}$
with $p$-symmetry.

We have
$$stab_{N_{Aut(F_n)}(B_k)}(\delta) = \left\{\matrix{
N_{\Sigma_p}(\Z/p) \times \Sigma_p \hfill
&\hbox{if } \delta \hbox{ has underlying graph }
 \hfill \cr
\hfill 
&\hfill \Theta_{p-1}^1  \hbox{or } \Theta_{p-1}^2 \cr
\hfill & \hfill \cr
N_{\Sigma_p}(\Z/p) \times H \hfill &\hbox{otherwise, where $H$ is a} 
\hfill \cr
\hfill &\hfill \hbox{group with } p \not | \hbox{ } |H| \cr} \right.$$

Arguments similar to those in Lemma \ref{t5} show that the $E_2$ 
page has
the form

$$E_2^{r,s} = \left\{\matrix{
\Z/p^2 \oplus \Z/p \hfill & r=0 \hbox{ and } s=0 \hfill \cr
H^r(\tilde Q_{p-1}; \Z/p) \hfill &1 \leq r \leq n-1, s = kn \hfill \cr
(2j-2)\Z/p \hfill & r=0 \hbox{ and } |s|=nj-1 \hfill \cr
(2j+1)\Z/p \hfill & r=0 \hbox{ and } |s|=nj \not = 0  \hfill \cr
0 \hfill &\hbox{otherwise } \hfill \cr} \right.$$
Let $j \geq 2$. Now the differential
$$d_{n} : E_{n}^{0,nj-1}
\to E_{n}^{n,n(j-1)}=0$$
is necessarily trivial and thus
$\hat H^{nj-1}(N_{Aut(F_n)}(B_{p-1}); \Z_{(p)})$ has a filtration
with successive terms 
$$E_2^{0,nj-1}=(2j-2)\Z/p$$
and
$$E_2^{2p-3,n(j-1)}=H^{2p-3}(\tilde Q_{p-1}; \Z/p).$$
Since
$$\hat H^{nj-1}(N_{Aut(F_n)}(B_{p-1}); \Z_{(p)})=
H^{nj-1}(N_{Aut(F_n)}(B_{p-1}); \Z_{(p)})$$
(because $nj-1$ is above the vcd of $N_{Aut(F_n)}(B_{p-1})$),
we can use the K\"unneth formula for the latter cohomology group to
specify the form that the above filtration takes and obtain that
$$\hat H^{nj-1}(N_{Aut(F_n)}(B_{p-1}); \Z_{(p)})
=(2j-2)\Z/p \oplus H^{2p-3}(\tilde Q_{p-1}; \Z/p).$$

The other tricky cohomology group to calculate is
$$\hat H^{-nj+1}(N_{Aut(F_n)}(B_{p-1}); \Z_{(p)})$$ (again,
for $j \geq 2$.)
This can be computed by noting that
$$d_2: E_2^{0,-nj+1} \to E_2^{2,-nj}$$
is zero by part 5 of Fact \ref{p1t1}
and that
$$E_2^{1,-nj}=H^1(\tilde Q_{p-1}; \Z/p)=0$$
by part 4 of Fact \ref{p1t1}.
The result follows for the case $k=p-1$. \END \end{pf*}

\begin{lemma} \label{t7}
$$\hat H^t(N_{Aut(F_n)}(C); \Z_{(p)}) \cong \left\{\matrix{
\Z/p^2 \oplus \Z/p \hfill &t = 0 \hfill \cr
({3k(p-1) \over 2} + 1) \Z/p \hfill &|t| = kn \not = 0 \hfill \cr
({3k(p-1) \over 2} - 2) \Z/p \hfill &|t| = kn-1 \hfill \cr
0 \hfill &\hbox{otherwise } \hfill \cr
} \right.$$
\end{lemma}

\begin{pf*}{Proof.} The normalizer $N_{Aut(F_n)}(C)$ acts on the 
contractible space
$X_n^C$ with finite stabilizers and finite quotient.
Hence we can use the equivariant cohomology
spectral sequence
(\ref{e2}) to calculate the cohomology of the normalizer.
This gives us:
\begin{equation} \label{e7}
E_1^{r,s} = \prod_{[\delta] \in \Delta^r}
\hat H^s(stab_{N_{Aut(F_n)}(C)}(\delta); \Z_{(p)}) \Rightarrow
\hat H^{r+s}(N_{Aut(F_n)}(C); \Z_{(p)})
\end{equation}
where $[\delta]$ ranges over the set $\Delta^r$ of orbits
of $r$-simplices $\delta$ in $X_n^C$.

\begin{claim} \label{tr37}
The quotient space 
$\bigcup_r\Delta^r$ has $3$ vertices and $2$ edges in it.  The 
vertices
correspond to marked graphs with underlying graphs $\Phi_{n}$,
$\Omega_{n}$, and $\Psi_{n}$.  The two edges come from the
forest collapses of $\Phi_{n}$ to $\Omega_{n}$ or
$\Psi_{n}$.  Pictorially, we have
$$\matrix{ & & \Phi_{n} & & \cr
           & \swarrow & & \searrow & \cr
           \Omega_{n} & & & & \Psi_{n} \cr}$$
\end{claim}

\begin{pf*}{Proof.} From Theorem \ref{tr4}, any two vertices 
(marked graphs)
of $X_n^C$ whose underlying graphs are reduced,
can be connected by a sequence of Nielsen
transformations.  The graphs that we obtain from
$\Psi_n$ by doing Nielsen moves are all isomorphic
to either $\Psi_n$ or $\Omega_n$.

It follows that if $\eta$ is a vertex of $X_n^C$ corresponding to
a reduced marked graph, then the underlying graph of
$\eta$ is either $\Psi_n$ or $\Omega_n$.

It remains to consider which graphs can be blowups
(see Definition \ref{tr22}) of $\Psi_n$ or
$\Omega_n$.  Such a blowup would have a nontrivial
$\Z/p$ action on at least $2p$ of its edges.  From this,
it is not hard to see (using similar methods
to those in Proposition \ref{t1}) that the only possibility
for the underlying graph of such a blowup is
$\Phi_n$. \END \end{pf*}
           
Direct examination reveals that the
vertex in $\Delta^0$ corresponding to $\Phi_{n}$ has
automorphism group $\Sigma_p \times \Z/2$.
In the notation used to define $\Phi_{n}$
(refer to the text just above Figure \ref{fig1}),
the $\Sigma_p$ in $\Sigma_p \times \Z/2$
acts on the collections of edges $\{a_i\}$,
$\{b_i\}$, and $\{c_i\}$, respectively,
by permuting their indices.  On the other hand,
the $\Z/2$ in $\Sigma_p \times \Z/2$
fixes the edges $a_i$ and switches the edges
$b_i$ with the edges $c_i$.
The group $C$ is included in $\Sigma_p \times \Z/2$
as the cyclic group generated by the
permutation $(1 2 \ldots p)$ in $\Sigma_p$.
Hence the subgroup of normalizing graph automorphisms in
$\Sigma_p \times \Z/2$ is
$$N_{\Sigma_p}(\Z/p) \times \Z/2.$$

The stabilizer of the vertex in $\Delta^0$
corresponding to $\Phi_n$ has cohomology
$$\hat H^*(N_{\Sigma_p}(\Z/p) \times \Z/2; \Z_{(p)})
=\hat H^*(\Sigma_p; \Z_{(p)}).$$

Similarly, the group of graph automorphisms of
$\Omega_{n}$ is $$\Sigma_p \times \Sigma_p.$$
The group $C$ is included in this as
the subgroup generated by
$$(1 2 \ldots p) \times (1 2 \ldots p).$$
The stabilizer of the vertex in $\Delta^0$ which
corresponds to $\Omega_n$ is 
the normalizer of $C$
in $\Sigma_p \times \Sigma_p$.
This normalizer is
$$(\Z/p \times \Z/p) \rtimes \Z/(p-1).$$
The generator of $\Z/(p-1)$ acts diagonally on
$(\Z/p \times \Z/p)$
by conjugating the generator of either $\Z/p$
to its $s$-th power for some generator $s$ of
$\F_p^\times$.

The cohomology of $(\Z/p \times \Z/p) \rtimes \Z/(p-1)$
can be calculated in a straightforward way from the cohomology of
$\Z/p \times \Z/p$ using Swan's theorem \ref{tr30} as
$$H^*((\Z/p \times \Z/p) \rtimes \Z/(p-1); \Z_{(p)}) =
H^*(\Z/p \times \Z/p; \Z_{(p)})^{\Z/(p-1)}.$$
For example, if $z_n, \bar z_n \in H^n(\Z/p; \Z/p)$ are the
generators corresponding to the first and second $\Z/p$'s in 
$\Z/p \times \Z/p$, respectively, then we can calculate the
cohomology of $(\Z/p \times \Z/p) \rtimes \Z/(p-1)$ in
dimensions $kn-1 > 0$ as follows.  For $i$ ranging from
$1$ to $k(p-1)-1$ elements of the form
$$z_{2i-1} \bar z_{2j}-z_{2i} \bar z_{2j-1}$$
are both $\Z/(p-1)$-invariant and are in the kernel of the
Bockstein homomorphism.  This gives $k(p-1)-1$ generators
for the cohomology group.
So we see that
$\hat H^t((\Z/p \times \Z/p) \rtimes \Z/(p-1); \Z_{(p)}) =$
$$\left\{\matrix{
\Z/p^2 \hfill &t = 0 \hfill \cr
(k(p-1) + 1) \Z/p \hfill &|t| = kn \not = 0 \hfill \cr
(k(p-1) - 1) \Z/p \hfill &|t| = kn-1 \hfill \cr
0 \hfill &\hbox{otherwise } \hfill \cr} \right.$$

Lastly, the
group of graph automorphisms of
$\Psi_{n}$ is
$$(\Sigma_p \times \Sigma_p) \rtimes \Z/2.$$
The group $C$ is included in this as
the subgroup generated by
$$(1 2 \ldots p) \times (1 2 \ldots p).$$
The stabilizer of the vertex in $\Delta^0$ which
corresponds to $\Psi_n$ is
the normalizer of $C$
in $(\Sigma_p \times \Sigma_p) \rtimes \Z/2.$
This normalizer is
$$(\Z/p \times \Z/p) \rtimes (\Z/2 \times \Z/(p-1)).$$
The generator of $\Z/(p-1)$ acts diagonally on
$(\Z/p \times \Z/p)$
by conjugating the generator of either $\Z/p$
to its $s$-th power for some generator $s$ of
$\F_p^\times$.  The $\Z/2$ acts by exchanging one
$\Z/p$ for the other in $\Z/p \times \Z/p$.
The cohomology of $(\Z/p \times \Z/p) \rtimes (\Z/2 \times \Z/(p-1))$
can be calculated using Swan's theorem \ref{tr30}
which indicates that
$$H^*((\Z/p \times \Z/p) \rtimes (\Z/2 \times \Z/(p-1)); \Z_{(p)}) =
H^*(\Z/p \times \Z/p; \Z_{(p)})^{(\Z/2 \times \Z/(p-1))}.$$
For example, we can calculate the
cohomology of $(\Z/p \times \Z/p) \rtimes (\Z/2 \times \Z/(p-1))$ in
dimensions $kn-1 > 0$ as follows.  For $i$ ranging from
$1$ to $k(p-1)/2-1$ elements of the form
$$z_{2i-1} \bar z_{2j}-z_{2i} \bar z_{2j-1} - z_{2j-1} \bar z_{2i}+
z_{2j} \bar z_{2i-1} $$
are $\Z/(p-1)$-invariant, $\Z/2$-invariant, and are in the 
kernel of the
Bockstein homomorphism.  This gives $k(p-1)/2-1$ generators
for the cohomology group.
Hence
$\hat H^t((\Z/p \times \Z/p) \rtimes (\Z/2 \times \Z/(p-1));
\Z_{(p)}) =$
$$\left\{\matrix{
\Z/p^2 \hfill &t = 0 \hfill \cr
({k(p-1) \over 2} + 1) \Z/p \hfill &|t| = kn \not = 0 \hfill \cr
({k(p-1) \over 2} - 1) \Z/p \hfill &|t| = kn-1 \hfill \cr
0 \hfill &\hbox{otherwise } \hfill \cr} \right.$$

The two edges in $\Delta^1$ have stabilizers isomorphic to
$N_{\Sigma_p}(\Z/p)$ or $N_{\Sigma_p}(\Z/p) \times \Z/2$.
We omit the argument here, but the stabilizers of
the edges can be found by examining which graph
automorphisms in
$$stab(\Phi_n) = N_{\Sigma_p}(\Z/p) \times \Z/2$$
preserve the relevant forest collapses.
In either case, if we take Farrell
cohomology with $\Z_{(p)}$-coefficients, then both edges have
stabilizers whose cohomology is the same as that of the symmetric
group $\Sigma_p$.

Combining all of this into the spectral sequence (\ref{e7}) and then
applying the differential on the $E_1$ page, we see that
$$E_2^{r,s} =
\left\{\matrix{
\Z/p^2 \oplus \Z/p \hfill &r = 0, s = 0\hfill \cr
({3k(p-1) \over 2} - 2) \Z/p \hfill &r=0, |s| = kn-1 \hfill \cr
({3k(p-1) \over 2} + 1) \Z/p \hfill &r=0, |s| = kn \not = 0 \hfill \cr
0 \hfill &\hbox{otherwise } \hfill \cr} \right.$$
Thus $E_2=E_\infty$ and $\hat H^*(N_{Aut(F_n)}(C); \Z_{(p)})$
is as stated. \END \end{pf*}

\begin{lemma} \label{t8}
$$\hat H^t(N_{Aut(F_n)}(D); \Z_{(p)}) \cong \hat H^t(\Sigma_p
\times \Sigma_p;\Z_{(p)}).$$
$$\hat H^t(N_{Aut(F_n)}(E); \Z_{(p)}) \cong
\hat H^t((\Sigma_p \times \Sigma_p) \rtimes \Z/2;\Z_{(p)})$$
$$ \cong \left\{\matrix{
\Z/p^2 \hfill &t = 0\hfill \cr
([k/2] + 1) \Z/p \hfill &|t| = kn \not = 0 \hfill \cr
([(k-1)/2]) \Z/p \hfill &|t| = kn-1 \hfill \cr
0 \hfill &\hbox{otherwise } \hfill \cr
} \right.$$
\end{lemma}

\begin{pf*}{Proof.} From Lemma \ref{tr11},
$$N_{Aut(F_n)}(D) = N_{\Sigma_p}(\Z/p) \times N_{\Sigma_p}(\Z/p).$$
and
$$N_{Aut(F_n)}(E) =
(N_{\Sigma_p}(\Z/p) \times N_{\Sigma_p}(\Z/p)) \rtimes \Z/2.$$

The detailed description of the cohomology of
$N_{Aut(F_n)}(E)$ in the statement of this lemma
is then obtained by using Swan's theorem \ref{tr30}. 
Generators in dimensions $kn-1 > 0$ 
come from expressions of the form
$$z_{in-1} \bar z_{jn}-z_{in} \bar z_{jn-1}$$
where $i$ ranges from $1$ to $[(k-1)/2]$.\END \end{pf*}

\section{Cohomology of $Aut(F_n)$: 
Proof of Theorem \ref{t9}} \label{c5}

In this section, we will use the lemmas of the previous section to
establish Theorem \ref{t9}.

\smallskip

\begin{pf*}{Proof of Theorem \ref{t9}.} 
The cohomology $\hat H^*(Aut(F_{n}); \Z_{(p)})$ will be
calculated using the normalizer spectral sequence (\ref{e3}), 
which has
$E_1$ page
$$E_1^{r,s} = \prod_{(A_0 \subset \cdots \subset A_r) \in 
|\mathcal{B}|_r}
\hat H^s( \bigcap_{i=0}^r N_{Aut(F_n)}(A_i); \Z_{(p)})
\Rightarrow \hat H^{r+s}(Aut(F_n); \Z_{(p)})$$
where $\mathcal{B}$ denotes the poset of conjugacy classes of
nontrivial elementary
abelian $p$-subgroups of $Aut(F_n)$, and $|\mathcal{B}|_r$ is the
set of $r$-simplices in $|\mathcal{B}|$.
We computed $|\mathcal{B}|$ in Proposition \ref{tr5}.  It is 
$1$-dimensional, so the above spectral sequence is zero except in
the columns $r=0$ and $r=1$.

Recall that the realization $|\mathcal{B}|$ of the poset 
$\mathcal{B}$ has $p$
path components.  One component just consists of a point corresponding
to the subgroup $A$.  In addition, $p-2$ other components are also
just points corresponding to the subgroups $B_k$ for
$k \in \{1,\ldots,p-2\}$.  Finally, the last component is a
$1$-dimensional simplicial complex corresponding to the subgroups
listed in diagram (\ref{e4}), which we duplicate here:
$$\matrix{ & & B_0 & & \cr
           & & \downarrow & & \cr
           & & D & & \cr
           & \nearrow & & \nwarrow & \cr
           B_{p-1} & & & & C \cr
           & \searrow & & \swarrow & \cr
           & & E & & \cr}$$

We have already calculated (in the lemmas of
the previous section) the contributions of all of
the vertices in $|\mathcal{B}|$ to the $E_1$ page in (\ref{e3}).

The contribution of a $1$-simplex in $|\mathcal{B}|$ can be obtained
by taking the cohomology of the intersections of the normalizers of
the vertices of the $1$-simplex.  Note that each of these
intersections is a finite group (since each is a subgroup of either
the finite group $N_{Aut(F_n)}(D)$ or $N_{Aut(F_n)}(E)$ of 
normalizing graph automorphisms.)
In this way, we can calculate the (now just Tate) cohomological
contributions of the $1$-simplices in (\ref{e4}) to be:
$$\matrix{ \hat H^*(N_{Aut(F_n)}(D) \cap N_{Aut(F_n)}(B_0); 
\Z_{(p)}) \hfill
&= \hat H^*(N_{\Sigma_p}(\Z/p)
 \times N_{\Sigma_p}(\Z/p);\Z_{(p)}) \hfill \cr
 \hfill &= \hat H^*(\Sigma_p \times \Sigma_p;\Z_{(p)}). \hfill}$$
$$\matrix{ \hat H^*(N_{Aut(F_n)}(D) \cap N_{Aut(F_n)}(B_{p-1}); 
\Z_{(p)}) \hfill
&= \hat H^*(N_{\Sigma_p}(\Z/p)
\times N_{\Sigma_p}(\Z/p);\Z_{(p)}) \hfill \cr
 \hfill &= \hat H^*(\Sigma_p \times \Sigma_p;\Z_{(p)}). \hfill}$$
$$\hat H^t(N_{Aut(F_n)}(D) \cap N_{Aut(F_n)}(C); \Z_{(p)}) =
\hat H^t((\Z/p \times \Z/p) \rtimes \Z/(p-1); \Z_{(p)}) \hfill$$
$$ =
\left\{\matrix{
\Z/p^2 \hfill &t = 0 \hfill \cr
(k(p-1) + 1) \Z/p \hfill &|t| = kn \not = 0 \hfill \cr
(k(p-1) - 1) \Z/p \hfill &|t| = kn-1 \hfill \cr
0 \hfill &\hbox{otherwise } \hfill \cr} \right.$$
$$\matrix{\hat H^*(N_{Aut(F_n)}(E) \cap N_{Aut(F_n)}(B_{p-1}); 
\Z_{(p)}) \hfill
&= \hat H^*(N_{\Sigma_p}(\Z/p)
\times N_{\Sigma_p}(\Z/p);\Z_{(p)}) \hfill \cr
 \hfill &= \hat H^*(\Sigma_p \times \Sigma_p;\Z_{(p)}).\hfill}$$
$$\hat H^t(N_{Aut(F_n)}(E) \cap N_{Aut(F_n)}(C); \Z_{(p)}) =
\hat H^t((\Z/p \times \Z/p) \rtimes (\Z/2 \times \Z/(p-1));
\Z_{(p)}) \hfill$$
$$ =
\left\{\matrix{
\Z/p^2 \hfill &t = 0 \hfill \cr
({k(p-1) \over 2} + 1) \Z/p \hfill &|t| = kn \not = 0 \hfill \cr
({k(p-1) \over 2} - 1) \Z/p \hfill &|t| = kn-1 \hfill \cr
0 \hfill &\hbox{otherwise } \hfill \cr} \right.$$

We are now ready to compute the $E_2$ page of the spectral sequence
(\ref{e3}).  The contributions coming from the isolated points of
$|\mathcal{B}|$ (i.e., from $A$,$B_1$, $\ldots$, $B_{p-2}$) survive
unaltered from the $E_1$ page.  The contributions from the 
$1$-dimensional
component
of $|\mathcal{B}|$ pictured in (\ref{e4}) will be what we concentrate
on from now on.  

First, we compute the values for the $E_2$ page in rows $s$
where $s=nj+k$ with $2 \leq k \leq 2(p-2).$  The
$E_1$ page is only nonzero in the column $r=0$ for these rows.  Hence
the entries in $E_1^{0,s}$ necessarily survive to the $E_2$ page and
from there survive to the $E_\infty$ page.  This gives us that 
$\hat H^t(Aut(F_{n}); \Z_{(p)})$ is as the proposition
claims for $t=nj+k$ with $2 \leq k \leq 2(p-2).$

For the rest of our calculations, we will use the fact that the
boundary map on the $E_1$ page is just the restriction map.  From a
comment by Brown in \cite{[B]} on page 286, we know that we can
compute these restriction maps (from normalizers of $p$-subgroups to
finite subgroups of those normalizers) by looking at the
$E_2$ pages of the various spectral sequences used to compute the
cohomologies of the normalizers (in the lemmas of the previous
section.)  This will help us to compute, for any row $s$, the value
$E_2^{0,s}$.  As an example of this, consider the copy
of $\hat H^*(N_{Aut(F_n)}(B_0); \Z_{(p)})$ contained in the column
$E_2^{0,s}$ of our spectral sequence (\ref{e3}).   Recall
from Lemma \ref{t6} that the cohomology of
$N_{Aut(F_n)}(B_0)$ is 
$$\hat H^t(N_{Aut(F_n)}(B_0); \Z_{(p)}) = 
\left\{\matrix{
\Z/p^2 \hfill &t = 0 \hfill \cr
(j+1) \Z/p \hfill &|t| = nj \not = 0 \hfill \cr
(j-1) \Z/p \hfill &|t| = nj-1 \hfill \cr
0 \hfill &|t| \equiv 1  \hbox{ } (\hbox{mod } n) \hfill \cr
0 \hfill &t \equiv 2  \hbox{ } (\hbox{mod } n) \hfill \cr
H^r(Q_{p-1};\Z/p) \hfill
&t \equiv r  \hbox{ } (\hbox{mod } n), \hfill \cr
 \hfill &3 \leq r \leq n-2 \hfill \cr
}\right.$$
and this was calculated by looking at a spectral
sequence whose $\mathcal{E}_2$ page
(which we now denote with a script
$E$ as $\mathcal{E}$ to distinguish it
from the spectral sequence (\ref{e3}) above)
is
$$\mathcal{E}_2^{r,s} = \left\{\matrix{
\Z/p^2 \hfill & r=0 \hbox{ and } s=0 \hfill \cr
H^r(Q_{p-1}; \Z/p) \hfill &1 \leq r \leq 2(p-2), s = kn \hfill \cr
(j-1)\Z/p \hfill & r=0 \hbox{ and } |s|=nj-1  \hfill \cr
(j+1)\Z/p \hfill & r=0 \hbox{ and } |s|=nj \not = 0  \hfill \cr
0 \hfill &\hbox{otherwise } \hfill \cr} \right.$$
From $(4.6)$ in Chapter X of \cite{[B]}, the restriction
map from $\hat H^*(N_{Aut(F_n)}(B_0); \Z_{(p)})$
to the ring $\frak{H}^*(N_{Aut(F_n)}(B_0); \Z_{(p)})$ 
deriving from the cohomology of the finite subgroups
of $N_{Aut(F_n)}(B_0)$ (see \cite{[B]} for a definition of this ring)
is the canonical surjection coming from the
vertical edge homomorphism
$$\hat H^t(N_{Aut(F_n)}(B_0); \Z_{(p)}) \to
\mathcal{E}_2^{0,t} =
\hat H^t(N_{\Sigma_p}(\Z/p) \times N_{\Sigma_p}(\Z/p);\Z_{(p)})$$
where
$$\hat H^t(N_{\Sigma_p}(\Z/p) \times N_{\Sigma_p}(\Z/p);\Z_{(p)})=
\left\{\matrix{
\Z/p^2 \hfill & t = 0 \hfill \cr
(j+1)\Z/p \hfill & |t| = nj \not = 0  \hfill \cr
(j-1)\Z/p \hfill & |t|= nj-1  \hfill \cr
0 \hfill &\hbox{otherwise } \hfill \cr} \right.$$
Recall that in order to know the column $E_2^{0,s}$ of
(\ref{e3}), we want to calculate the
restriction map
$$\matrix{
\hat H^*(N_{Aut(F_n)}(B_0); \Z_{(p)})  \to 
\hat H^*(N_{Aut(F_n)}(B_0) \cap N_{Aut(F_n)}(D)) \hfill \cr
\hfill
\hat H^*(N_{\Sigma_p}(\Z/p) \times N_{\Sigma_p}(\Z/p);\Z_{(p)}) \cr}$$
This restriction map is
the vertical edge homomorphism calculated above.

We
can similarly
look at the $\mathcal{E}_2$-pages of the other spectral sequences in
section \ref{c4} to calculate the
restriction maps from the cohomologies of
$N_{Aut(F_n)}(B_{p-1})$ and $N_{Aut(F_n)}(C)$ to their finite
subgroups
$$N_{Aut(F_n)}(B_{p-1}) \cap N_{Aut(F_n)}(D), \hbox{ }
N_{Aut(F_n)}(B_{p-1}) \cap N_{Aut(F_n)}(E),$$
$$N_{Aut(F_n)}(C) \cap N_{Aut(F_n)}(D), \hbox{ and }
N_{Aut(F_n)}(C) \cap N_{Aut(F_n)}(E).$$
In these two cases, the restriction map is not just the vertical
edge homomorphism.  This is because the $E_1$-pages
(and hence the $E_2$-pages) of the spectral
sequences used to calculate 
$N_{Aut(F_n)}(B_{p-1})$ and $N_{Aut(F_n)}(C)$ each had vertical 
edges with
cohomology groups coming from two different underlying marked
graphs, namely $\Omega_{2(p-1)}$ and $\Psi_{2(p-1)}$.  The
cohomological contribution from  $\Omega_{2(p-1)}$ is what we are 
concerned 
with when calculating the restriction to 
$N_{Aut(F_n)}(B_{p-1}) \cap N_{Aut(F_n)}(D)$
or $N_{Aut(F_n)}(C) \cap N_{Aut(F_n)}(D)$, while the
contribution from  $\Psi_{2(p-1)}$ is what we are concerned with
when calculating restrictions to intersections involving
$N_{Aut(F_n)}(E)$.  So the restriction maps in these cases
are found by composing the vertical edge homomorphism with a 
restriction
map onto the portions of the $E_2$ page that come from either
$\Omega_{2(p-1)}$ or $\Psi_{2(p-1)}$.
On the other hand, the restrictions from
$N_{Aut(F_n)}(D)$ and $N_{Aut(F_n)}(E)$ to their subgroups are easy to
compute as $N_{Aut(F_n)}(D)$ and $N_{Aut(F_n)}(E)$ are just well
known finite groups.

After calculating all of
the terms in the column $E_2^{0,s}$ of (\ref{e3})
as  above, the value $E_2^{1,s}$ 
is found from 
an Euler characteristic argument using $E_1^{0,s}$ and
$E_1^{1,s}$.  That is, if $s \not = 0$ then $E_1^{0,s}$ and
$E_1^{1,s}$ are the only nonzero terms on the row $s$ and they
are both $\F_p$-vector spaces.  As the boundary map goes from
$E_1^{0,s} \to E_1^{1,s}$, this yields
$dim_{\F_p}(E_1^{0,s}) - dim_{\F_p}(E_1^{1,s})
= dim_{\F_p}(E_2^{0,s}) - dim_{\F_p}(E_2^{1,s})$
by the standard Euler characteristic argument.

Assume $s=-kn+1<0$ (the case $s=kn-1>0$
follows similarly, with the only exception being the additional
summand of $H^{n-1}(\tilde Q_{p-1}; \Z/p)$ that needs
to be dealt with.)  Then
$E_2^{0,s} = (k+[(k-1)/2]-1)\Z/p$.  Now $k-1$ of these $\Z/p$'s 
come from
$0$-cocycles that are summations of cocycles from:
\begin{itemize}
\item The portion of $\hat H^s(N_{Aut(F_n)}(B_{p-1});\Z_{(p)})$ 
coming from the
graph $\Omega_n$. (In the spectral sequence we used to calculate
$\hat H^s(N_{Aut(F_n)}(B_{p-1});\Z_{(p)})$ in part $3.$ of
Lemma \ref{t6}, this was the contribution given by the
stabilizer of the graph $\Theta_{p-1}^2$.)
\item $\hat H^s(N_{Aut(F_n)}(D);\Z_{(p)}).$
\item $\hat H^s(N_{Aut(F_n)}(B_{0});\Z_{(p)}).$
\item The portion of $\hat H^s(N_{Aut(F_n)}(C);\Z_{(p)})$ coming 
from the
graph $\Omega_n$ (in the spectral sequence we used to calculate
$\hat H^s(N_{Aut(F_n)}(C));\Z_{(p)})$.)
\end{itemize}
The other $[(k-1)/2]$ of the $\Z/p$'s come from $0$-cocycles that are
summations of cocycles from:
\begin{itemize}
\item The portion of $\hat H^s(N_{Aut(F_n)}(B_{p-1});\Z_{(p)})$ 
coming from the
graph $\Psi_n$. 
(In the spectral sequence we used to calculate
$\hat H^s(N_{Aut(F_n)}(B_{p-1});\Z_{(p)})$ in part $3.$ of
Lemma \ref{t6}, this was the contribution given by the
stabilizer of the graph $\Theta_{p-1}^1$.)
\item $\hat H^s(N_{Aut(F_n)}(E);\Z_{(p)}).$
\item The portion of $\hat H^s(N_{Aut(F_n)}(C);\Z_{(p)})$ coming 
from the
graph $\Psi_n$ (in the spectral sequence we used to calculate
$\hat H^s(N_{Aut(F_n)}(C);\Z_{(p)})$.)
\end{itemize}
This allows us to find that the contribution of the component of
$|\mathcal{B}|$
corresponding to diagram (\ref{e4}) in the row $s=-kn+1<0$
is
$E_2^{0,s} = (k+[(k-1)/2]-1)\Z/p.$
Since
$dim_{\F_p}(E_1^{0,s})-dim_{\F_p}(E_1^{1,s})=$
$$(4(k-1)+{3k(p-1) \over 2}+\left[{k-1 \over 2}\right]-2) -
(3(k-1)+{3k(p-1) \over 2}-1)=k+\left[{k-1 \over 2}\right]-1,$$
we find that $E_2^{1,s}=0$.

Assume $s=kn>0$ or $s=-kn<0$.  Then
$E_2^{0,s} = (k+[k/2]+p-1)\Z/p$.  First, observe that $p-1$ of these
$\Z/p$'s come from the points in $|\mathcal{B}|$ corresponding to
$A$,$B_1$, $\ldots$, $B_{p-2}$.
Also, $k-1$ of these $\Z/p$'s come from summations of
cocycles in
$H^s(N;\Z_{(p)})$ for $N$ equal to
$N_{Aut(F_n)}(B_0)$, $N_{Aut(F_n)}(B_{p-1})$, $N_{Aut(F_n)}(C)$, and
$N_{Aut(F_n)}(E)$ where in all cases the cohomology
came from the graph $\Omega_n$ (in the spectral sequences
used to calculate the cohomology groups of the various
normalizers.)
In addition, $[k/2]$ of the $\Z/p$'s come from $0$-cocycles that are
summations of cocycles from $H^s(N;\Z_{(p)})$ for $N$ equal to
$N_{Aut(F_n)}(B_{p-1})$, $N_{Aut(F_n)}(C)$, and $N_{Aut(F_n)}(E)$ 
where in all cases the cohomology
came from the graph $\Psi_n$.
Finally, one cocycle comes from summing up cocycles from
$H^s(N;\Z_{(p)})$ for $N$ equal to $N_{Aut(F_n)}(B_0)$, 
$N_{Aut(F_n)}(B_{p-1})$,
$N_{Aut(F_n)}(C)$, $N_{Aut(F_n)}(D)$, and $N_{Aut(F_n)}(E)$.
This allows us to find that the contribution of the component of $B$
corresponding to diagram (\ref{e4}) in the row $s=kn>0$ or $s=-kn<0$
is
$E_2^{0,s} = (k+[k/2]+p-1)\Z/p.$
Note that because
$dim_{\F_p}(E_1^{0,s})-dim_{\F_p}(E_1^{1,s})=$
$$(4k+{(3k+2)(p-1) \over 2}+[k/2]+5) -
(3k+{3k(p-1) \over 2}+5)=k+[k/2]+p-1,$$
we find that $E_2^{1,s}=0$.

Next, $E_2^{0,0}$ of the normalizer spectral sequence 
is readily computed to be $\Z/p^2 \oplus p(\Z/p)$,
where $p-1$ of the $\Z/p$'s come from the isolated points of
$|\mathcal{B}|$ and $\Z/p^2 \oplus \Z/p$ comes from vertices in
the $1$-dimensional component of $|\mathcal{B}|$.
In addition, $E_2^{0,1} = \Z/p$ where
the $\Z/p$ corresponds to edges of $|\mathcal{B}|$ which give
cohomology classes in $E_1^{0,1}$ that are not mapped
onto by cohomology classes in $E_1^{0,0}$.

We further illustrate the above calculations by looking at
the contribution of diagram (\ref{e4}) to the explicitly listed
(i.e., not listed at coming from cohomologies of various quotient
spaces) cohomology classes in $E_2^{0,*}$.  This the same as 
that contributed by the
by the graph of groups listed in the figure below.

\bigskip

\centerline{\input{cohom3.pic}}
\begin{figure}[here]
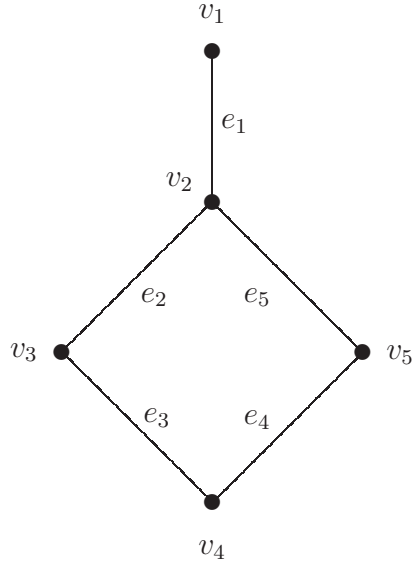

\caption{\label{gogps} A graph of groups}
\end{figure}

In the graph of groups in Figure \ref{gogps}, there are
five vertices $v_i$ and five edges $e_i$. The vertex 
and edge groups are
\begin{itemize}
\item $v_1$: $N_{\Sigma_p}(\Z/p) \times \Sigma_p$
\item $v_2$: $N_{\Sigma_p}(\Z/p) \times N_{\Sigma_p}(\Z/p)$
\item $v_3$: $(N_{\Sigma_p}(\Z/p) \times \Sigma_p) \ast_{\Sigma_p}
(N_{\Sigma_p}(\Z/p) \times \Sigma_p)$, where the amalgamated 
$\Sigma_p$
includes into the $\Sigma_p$ on the right in each factor.
\item $v_4$:  $(N_{\Sigma_p}(\Z/p) \times N_{\Sigma_p}(\Z/p)) 
\rtimes \Z/2$
\item $v_5$: $\left[ (\Z/p \times \Z/p) \rtimes \Z/(p-1) \right] 
\ast_{N_{\Sigma_p}(\Z/p)}
\left[ (\Z/p \times \Z/p) \rtimes (\Z/2 \times \Z/(p-1)) \right]$,
where the amalgamated $N_{\Sigma_p}(\Z/p)$ includes diagonally into
the factors on both sides.
\item $e_1$: $N_{\Sigma_p}(\Z/p) \times N_{\Sigma_p}(\Z/p)$
\item $e_2$: $N_{\Sigma_p}(\Z/p) \times N_{\Sigma_p}(\Z/p)$
\item $e_3$: $ N_{\Sigma_p}(\Z/p) \times N_{\Sigma_p}(\Z/p)$
\item $e_4$: $ (\Z/p \times \Z/p) \rtimes (\Z/2 \times \Z/(p-1))$
\item $e_5$: $ (\Z/p \times \Z/p) \rtimes \Z/(p-1)$
\end{itemize}
The vertex groups consist of the explicitly listed
contributions from the normalizers of
$B_0$, $D$, $B_{p-1}$, $E$, and $C$, respectively, 
to the column $E_2^{0,*}$ of the normalizer spectral sequence.
For example, from the action of $N_{\Sigma_p}(\Z/p) \times \Sigma_p$
on the graph $\Omega_n$, we obtain an inclusion of this group into
$N_{Aut(F_n)}(B_0)$.  We similarly obtain two inclusions of
$N_{\Sigma_p}(\Z/p) \times \Sigma_p)$ into $N_{Aut(F_n)}(B_{p-1})$
by taking actions on $\Omega_n$ or $\Psi_n$, respectively.  This
yields (cf Case 3 of the proof of Lemma \ref{t6}) an inclusion of
the vertex group for $v_3$ into $N_{Aut(F_n)}(B_{p-1})$. 
Lastly, from the proof of Lemma \ref{t7}, we see that 
$N_{Aut(F_n)}(C)$ is isomorphic to the vertex group for $e_5$.
The edge groups above  are just the intersections of the 
normalizers of the $B_0$, $D$, $B_{p-1}$, $E$, and $C$. 

Using the property that $X \ast_{Y} Y \cong X$ and noting that
the cohomology of $\Sigma_p$ and $N_{\Sigma_p}(\Z/p)$ are the same,
we see that graph of groups in Figure \ref{gogps} can be reduced 
to the graph of groups listed in Figure \ref{gogps2} which yielded
the group $\varsigma_p$. 
The left hand vertex in the graph of groups in Figure \ref{gogps2}
comes from the first half of the amalgamated product 
decomposition of the vertex group $v_3$, and 
the right hand vertex comes from $v_4$.  The top edge group comes
from the $\Sigma_p$ being amalgamated over to make $v_3$, and the
bottom edge group comes from the $N_{\Sigma_p}(\Z/p)$ being 
amalgamated over to make $v_5$.

For example, we can look at the contribution of 
Figure \ref{gogps2} to the row $s=n$ by examining an equivariant 
spectral
sequence for $\varsigma_p$.  This would have 3 generators
in $E_1^{0,n}$, two from $\Sigma_p \times \Sigma_p$ and one 
from $(\Sigma_p \times \Sigma_p) \rtimes \Z/2$.  It would also
have 2 generators in $E_1^{1,n}$, one each from the top edge group 
in Figure \ref{gogps2} and one from the bottom edge group.
The 2 by 3 matrix corresponding to the coboundary map in this row
has rank 2, since the first generator of $E_1^{0,n}$ includes only
into the second generator of $E_1^{1,n}$ while the second generator
of $E_1^{0,n}$ includes diagonally.  Since the matrix
has a $1$-dimensional nullspace, diagram (\ref{e4}) contributes
one $\Z/p$ to $E_2^{0,n}$ and nothing to $E_2^{1,n}$ is $0$.
Other rows can be examined similarly.

Finally, the last remark in the statement of the theorem 
is justified since all of the
explicitly listed cohomology classes in the column $E_2^{0,*}$
coming from isolated points of $|\mathcal{B}|$ 
are detected upon restriction to the
portions of the cohomologies of the normalizers in $Aut(F_n)$ 
of $A$, $B_1$,
$\ldots$, $B_{p-2}$, coming from the stabilizers of
marked graphs corresponding to $R_{2(p-1)}$, 
$R_1 \vee \theta_{p-1} \vee R_{p-2}$, $\ldots$, 
$R_{p-2} \vee \theta_{p-1} \vee R_{1}$.  The 
explicitly listed classes in $E_2^{0,*}$ coming
from diagram (\ref{e4}) all come from the graph of groups 
$\varsigma_p$ listed
in Figure \ref{gogps2}.  In the figure, the 
$\Sigma_p \times \Sigma_p$ comes from graph automorphisms
of $\Omega_{2(p-1)}$ and the 
$(\Sigma_p \times \Sigma_p) \rtimes \Z/2$ comes from graph 
automorphisms of 
$\Psi_{2(p-1)}$. 
The fact that the
isomorphism from the ``explicit'' part of the cohomology
of $Aut(F_n)$ to that of
$(\raisebox{-1.1ex}{\biggie *}\hbox{}_{i=1}^{p-1} \Sigma_p) 
\ast \varsigma_p$ preserves the ring structure follows from
(4.5) (vi) in Chapter X of \cite{[B]}. \END \end{pf*}

\appendix

\section{The usual cohomology of $Aut(F_{2(p-1)}$)} \label{c6}

We now state the (only partial) 
characterization of $H^t(Aut(F_{n}); \Z_{(p)})$ 
which results from our calculations in Theorem \ref{t9},
where $n=2(p-1)$ as before.

\begin{prop} \label{t10}
$$\matrix{
\hfill H^t(Aut(F_{n}); \Z_{(p)}) &=& &\left\{\matrix{
\hat H^t(Aut(F_{n}); \Z_{(p)}) \hfill &t > 4p-6 \hfill \cr
H^t(Q_{n}; \Z_{(p)}) \hfill &0 \leq t \leq 2p-3 \hfill \cr}
\right. \hfill \cr }$$
\end{prop} 

\begin{pf*}{Sketch of proof.}
Since the virtual cohomological dimension of
$Aut(F_{n})$ is $4p-6$, it follows directly (see \cite{[B]}) that
$H^t(Aut(F_{n}); \Z_{(p)}) = \hat H^t(Aut(F_{n});
\Z_{(p)})$
for $t > 4p-6$ and that the sequence
$$H^{m}(\Gamma); \Z_{(p)}) \to H^{m}(Aut(F_{n});
\Z_{(p)})
\to \hat H^{m}(Aut(F_{n}); \Z_{(p)}) \to 0$$
is exact for $m = 4p-6$, where $\Gamma$ is any torsion free subgroup
of $Aut(F_{n})$ of finite index
and
where the first map is the transfer map.

To establish the final part of the proposition,
we can use the equivariant cohomology 
spectral sequence (\ref{e1}) to calculate
$H^t(Aut(F_{n}); \Z_{(p)})$ for $0 \leq t \leq 2p-3$.
We want to show that for $0 < s \leq 2p-3$,
the rows $\tilde E_2^{t,s}$ are all zero, so that 
$H^s(Aut(F_{n}); \Z_{(p)}) = H^s(Q_{n}; \Z_{(p)})$.    

Let $X_s$ be
the $p$-singular locus of $X_{n}$, or the set of all simplices
in $X_{n}$
such that $p$ divides $|stab(\delta)|$. Let $\frak{E}_1^{r,s}$
be the $E_1$-page of the equivariant cohomology spectral sequence
for calculating $H_{Aut(F_{n})}^*(X_s;\Z_{(p)})$. 
Note that for $s > 0$, $\tilde E_1^{*,s} = \frak{E}_1^{r,s}$
and $\tilde E_2^{*,s} = \frak{E}_2^{r,s}$
Further note that $X_s$ separates into $p$ disjoint,
$Aut(F_{n})$-invariant (not necessarily connected)
subcomplexes corresponding to the $p$ path components
of $|\mathcal{B}|$.  (These path components were mentioned in
the proof of Theorem \ref{t9}.)
So we see that
$$H_{Aut(F_{n})}^*(X_s; \Z_{(p)})
= \prod_Y \hat H_{Aut(F_{n})}^*(Y; \Z_{(p)})$$
where $Y$ ranges over the $p$ distinct, disjoint,
$Aut(F_{n})$-invariant subcomplexes.  Hence the $E_1$
and $E_2$
pages can be calculated separately for each $Y$.  
Only one of the subcomplexes $Y$ has dihedral
symmetry in it and is relevant to the horizontal rows
between $1$ and $2p-3$;  
namely, the subcomplex $Y_A$ which has
marked graphs with underlying graph the
rose $R_{n}$ in it. Observe that $Y_A$ is
$Aut(F_{n}) \cdot X_{n}^A$, several disjoint copies of
$X_{n}^A$ grouped together.  Let $\Xi_p$ be the
$1$-skeleton of the cone over a
$p$-gon.  Note that $\Xi_p$ has dihedral symmetry. 
It is easily shown that the only marked graphs 
with dihedral symmetry are in $Y_A$ and have underlying graph
$\Xi_p \vee \Gamma_{p-2}$,  where $\Gamma_{p-2}$ is any
pointed graph with fundamental group isomorphic to $F_{p-2}$, and
where the
wedge does not necessarily take place at the basepoint.

Now consider a specific path component of $Y_A$, say $X_{n}^A$.
There is an action of $A \cong \Z/p$ on each simplex of
$X_{n}^A$.  In \cite{[K-V]}, Krstic and Vogtmann define an
equivariant deformation retract $L_A$ of $X_{n}^A$ by only
keeping
the ``essential graphs'' in $X_{n}^A$.  Among the graphs that are
in $X_{n}^A - L_A$ are all of
the ones with dihedral symmetry, along with various other
inessential graphs.
Define $L_A'$ to be the subcomplex of
$X_{n}^A$ obtained by collapsing the wedge summand
$\Xi_p$ to $R_p$ in graphs with dihedral symmetry.  
By the same poset lemma used in \cite{[K-V]},
the collapse from $X_{n}^A$
to $L_A'$ is an equivariant deformation retraction.
Observe that
$X_{n}^A \supset L_A' \supset L_A$.

Denote the $E_1$-page of
the equivariant cohomology spectral sequence used to calculate
$\hat H_{N_{Aut(F_n)}(A)}^*(X_{n}^A; \Z_{(p)})$
by  $\mathcal{E}_1^{r,s}$.
The contribution of the equivariant cohomology of $Y_A$
to $\frak{E}_1^{r,s}$ is the same as $\mathcal{E}_1^{r,s}$.
For $0 < s \leq 2p-3$ and
$s \not = 4k$, the row $\mathcal{E}_1^{*,s}$ is all zero.
On the other hand, for $0 < s=4k < 2p-3$, direct examination
reveals that the
row $\mathcal{E}_1^{*,s}$ is the relative cochain
complex $C^*(X_{n}/N_{Aut(F_n)}(A),L_A'/N_{Aut(F_n)}(A);\Z/p)$.
Because the homotopy
in the deformation retraction from $X_{n}^A$ to $L_A'$
is $N_{Aut(F_n)}(A)$-equivariant,
the
relative cohomology groups
$$H^*(X_{n}/N_{Aut(F_n)}(A),L_A'/N_{Aut(F_n)}(A);\Z/p)$$
are all zero.  Hence all of the rows $\mathcal{E}_2^{r,s}$
are zero for $0 < s \leq 2p-3$, which completes
the proof. \END \end{pf*}

\section{The Farrell cohomology of $Aut(F_l)$ \newline
for $l < 2(p-1)$} \label{c7}

For $l < p-1$, $\hat H^*(Aut(F_l); \Z_{(p)}) = 0$ from  an 
easy spectral sequence
argument.
We have already remarked that Glover and Mislin's work
in \cite{[G-M]} directly implies
that 
$$\hat H^*(Aut(F_{p-1}); \Z_{(p)}) = \hat H^*(\Sigma_p; \Z_{(p)}),$$
$$\hat H^*(Aut(F_{p}); \Z_{(p)}) = 3\hat H^*(\Sigma_p; \Z_{(p)}),$$
and we have noted that Yu Qing Chen's work in \cite{[C]}
shows that 
$$\hat H^*(Aut(F_{p+1}); \Z_{(p)}) = 4\hat H^*(\Sigma_p; \Z_{(p)})$$
and that 
$\hat H^t(Aut(F_{p+2}); \Z_{(p)}) = 5\hat H^t(\Sigma_p; \Z_{(p)})
\oplus \hat H^{t-4}(\Sigma_p; \Z_{(p)})$.

In this section, we show how to 
calculate $\hat H^*(Aut(F_{l}); \Z_{(p)})$
for $p \leq l < n$, where $n=2(p-1)$.
We do this by modifying in a direct manner
the arguments we have made
in the previous sections to calculate
$\hat H^*(Aut(F_{2(p-1)}); \Z_{(p)})$.

As before, one of the elements of $\mathcal{B}$
corresponds to the subgroup $A \cong \Z/p$ of
$Aut(F_l)$ obtained by rotating the first
$p$ petals of the rose $R_l$.
The other elements of $\mathcal{B}$
correspond to subgroups $B_k \cong \Z/p$
for $k \in \{0, \ldots, l-p+1\}$.
The subgroup $B_k$ is obtained,
as before, by rotating the $p$-edges
of the $\Theta$-graph in the middle of
$R_k \vee \Theta_{p-1} \vee R_{l-p+1-k}$.

As in Lemma \ref{t5}, we see that
$$N_{Aut(F_l)}(A) \cong
N_{\Sigma_p}(\Z/p) \times
((F_{l-p} \times F_{l-p}) \rtimes (\Z/2 \times Aut(F_{l-p})).$$
Let $\langle \omega \rangle  \cong \Z/2$ be the subgroup of
$Aut(F_{l-p+2})$ corresponding to the
action given by
switching the first two petals of the rose $R_{l-p+2}$.
Accordingly,
$$N_{Aut(F_l)}(A) \cong C_{Aut(F_l)}(A) \cong
(F_{l-p} \times F_{l-p}) \rtimes (\Z/2 \times Aut(F_{l-p}).$$
As in Definition \ref{tr17}, let
$X_{l-p+2}^{\langle \omega \rangle}$ be the fixed point set of 
$\omega$
in $X_{l-p+2}$.  Then let $X_{l-p+2}^\omega$
be the deformation retract of $X_{l-p+2}^{\langle \omega \rangle}$
obtained by collapsing out inessential edges of
marked graphs.  Finally, define $Q_{l-p+2}^\omega$
to be the quotient of $X_{l-p+2}^\omega$
by $N_{Aut(F_{l-p+2})}(\omega)$.

Putting all of this together, and using the same
methods as those used in Lemma \ref{t5},
Lemma \ref{t6} and Theorem \ref{t9}, we see that
$$\hat H^t(N_{Aut(F_l)}(A); \Z_{(p)})
= \left\{\matrix{
\Z/p \hfill &t \equiv 0  \hbox{ } (\hbox{mod } n) \hfill \cr
H^r(Q_{l-p+2}^\omega; \Z/p) \hfill &t \equiv r  \hbox{ } 
(\hbox{mod } n), \hfill \cr
 &2 \leq r \leq 2(l-p) \hfill \cr
0 \hfill &\hbox{otherwise } \hfill \cr} \right.$$
and that for $k \in \{0,\ldots,l-p+1\}$,
$\hat H^t(N_{Aut(F_l)}(B_k); \Z_{(p)}) = $
$$\matrix{
\hfill \cr
\left\{\matrix{
\Z/p \hfill &t \equiv 0  \hbox{ } (\hbox{mod } n)  \hfill \cr
H^r(\tilde Q_{k} \times Q_{p-1-k};\Z/p) \hfill
&t \equiv r  \hbox{ } (\hbox{mod } n), \hfill \cr
\hfill &2 \leq r \leq 2(l-p)+1 \hfill \cr
0 \hfill &\hbox{otherwise } \hfill \cr} \right.
\hfill \cr
\hfill & \hfill \cr}$$
so we have that for $p$ odd and $l \in \{p, \ldots, 2p-3\}$,
$\hat H^t(Aut(F_{l}); \Z_{(p)}) = $
$$\matrix{
\hfill \cr
\left\{\matrix{
(l-p+3)\Z/p \hfill &t \equiv 0  \hbox{ } (\hbox{mod } n)  \hfill \cr
\hfill & \hfill \cr
\sum_{k=0}^{l-p+1} H^r(\tilde Q_{k} \times Q_{l-p+1-k};\Z/p)
\hfill &t \equiv r  \hbox{ } (\hbox{mod } n), \hfill \cr
\oplus H^r(Q_{l-p+2}^\omega;\Z/p) \hfill &2 \leq r \leq 2(l-p)+1 
\hfill \cr
\hfill & \hfill \cr
0 \hfill &\hbox{otherwise } \hfill \cr} \right. \hfill \cr}$$

\smallskip

The aforementioned results of Glover and Mislin
and Chen give more information than our results
here, however, as they actually explicitly calculate the
cohomology groups
$$\sum_{k=0}^{l-p+1} H^r(\tilde Q_{k} \times Q_{l-p+1-k};\Z/p)
\oplus H^r(Q_{l-p+2}^\omega;\Z/p)$$
that arise in the above formula in their cases.
        
\section{On the prime $p=3$: Proof of Corollary \ref{tapp}} 
\label{appen}

For the sake of having a concrete example, we
calculate the cohomologies of all of the
quotient spaces involved in Theorem \ref{t9} when $p=3$.
This was also done, independently, by Glover and
Henn.

\begin{pf*}{Proof of Corollary \ref{tapp}.}
Examining Theorem \ref{t9} reveals that we
must show that
none of the various groups
$H^r(\tilde Q_{k}; \Z/3)$, $H^r(Q_{k}; \Z/3)$,
and $H^r(Q_{3}^\omega; \Z/3)$ (where $k = 1,2$)
contribute any nonzero cohomology classes.

The groups $H^r(Q_{1}; \Z/3)$, $H^r(Q_{2}; \Z/3)$,
$H^r(\tilde Q_{1}; \Z/3)$, and $H^r(\tilde Q_{2}; \Z/3)$
are all zero by Fact \ref{p1t1}.

For $H^r(Q_{3}^\omega; \Z/3)$
(see Definition \ref{tr17}
to recall the definitions related to this space), 
we note that the relevant
marked graphs in 
$Q_{3}^{\langle \omega \rangle} \cong X_4^A / N_{Aut(F_4)}(A)$
are those listed in component $(\mathcal{A})$
of section 4 in the paper \cite{[G-M]} by Glover and Mislin,
with the additional complication that a basepoint $*$ can
be added to the graphs in various places.  However, most of these
graphs have inessential edges (see \cite{[K-V]}) under the action
of $N_{Aut(F_4)}(A)$, and are thus collapsed directly
away when we reduce from  
$Q_{3}^{\langle \omega \rangle} \cong X_4^A / N_{Aut(F_4)}(A)$
to the space $Q_{3}^\omega$.
We list the graphs
from \cite{[G-M]} that give
$X_4^A / N_{Aut(F_4)}(A)$
here:
\begin{itemize}
\item $R_4$.  The rose has no inessential edges, even when you
attach the basepoint $*$ to the middle of one of the petals.
\item $\Theta_4$.  This $\Theta$-graph has no inessential
edges, regardless of where the basepoint is attached.
\item $W_3 \vee R_1$.  In our notation, this would be
$\Xi_3 \vee R_1$. The inessential edges in the ``spokes'' of the
graph $W_3$ are collapsed and reduce this graph to $R_4$.
\item $\Theta_3 * R_1$.  This is $\Theta_3$ with a loop $R_1$
attached to the middle of one of the edges of the $\Theta$-graph.
The edge of the $\Theta$-graph that the loop is attached
to is inessential.  (It will be 2 or 3 actual edges in the
resulting graph, all of which are inessential, depending upon
where the basepoint is placed.)  Collapsing the
inessential edges yields $R_4$.
\item $\Theta_2 \lozenge Y$.  A graph in the shape of a
letter $Y$ attached to the graph $\Theta_2$, with the
top vertices of the $Y$ attached to one side of the
$\Theta$-graph, and the bottom vertex of the $Y$ attached
to the other side of the $\Theta$-graph.  The bottom edge of
the $Y$ is inessential.  Collapsing this gives $\Theta_4$.
\item $\Theta_2 ** \Theta_1$.  Two $\Theta$-graphs with
a line drawn from the left vertex of one to the left
vertex of the other, and a line drawn from the right vertex of
one to the right vertex of the other.  The new lines drawn
are inessential edges, and can be collapsed away to
yield $\Theta_4$.
\end{itemize}

So we are left with 4 basepointed graphs.  Two come
from the rose $R_4$, depending upon where we place the
basepoint, and the other two come from $\Theta_4$ in
a similar manner.  In particular, only $\Theta_4$
(with the basepoint $*$ placed in the middle of one of
its edges)
can contribute a $2$-simplex to our complex, and the
relevant marked graph only has one maximal subforest
(up to an isomorphism of the graph).  Hence it contributes
exactly two $2$-simplices, which join together
to form a square.
Consequently it is clear that $H^2(Q_{3}^\omega; \Z/3)=0$,
which is all we needed to show to prove that
$Q_{3}^\omega$ contributes nothing more to our
cohomology calculations.
  
\bigskip

\input{cohom4.pic}
\begin{figure}[here]
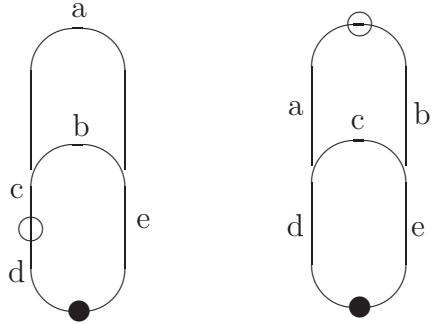

\caption{\label{fg6} Graphs giving some $3$-simplices}
\end{figure}

For the final case of
considering the contributions of $H^3(\tilde Q_{2}; \Z/3)$,
we again use arguments like
those in \cite{[V]} and Proposition $10.3$ of \cite{[JD]}.
We show that all of the $3$-simplices in 
$\tilde Q_{2}$ can be collapsed away,
so that $H^3(\tilde Q_{2}; \Z/3)$ is
necessarily zero.  The relevant graphs which can
give $3$-simplices are listed in Figure \ref{fg6}.
For each graph, the filled-in dot is the basepoint, but 
the open circle $\circ$
is the other ``distinguished point'' of the
graph 
and indicates where a $\Theta$-graph $\Theta_2$
should be attached.

The first of these graphs
has four subforests $\{b,c,e\}$, $\{b,c,d\}$,
$\{b,d,e\}$, and $\{c,d,e\}$,
each of which gives a
$3$-dimensional cube.
These can be collapsed in a manner
similar to that described
in Proposition $10.3$ in \cite{[JD]}.
That is, the cube corresponding to the
first subforest $\{b,c,e\}$
has a free minusface obtained by collapsing $b$.
So we can collapse the interior of this cube away from that face.
Then the cube $\{b,c,d\}$ has a free
plusface corresponding to
$\{b,c\}$ and $\{b,d,e\}$ has a free
plusface corresponding to
$\{b,e\}$.  Both of these cubes can
be collapsed away from those
respective plusfaces.
This leaves the cube corresponding to
$\{c,d,e\}$ with all plusfaces free.
Thus we can disregard the first of the graphs
in Figure \ref{fg6}.
(For a more detailed
description of what plusfaces and minusfaces are,
along with several more examples, see Chapter 10 of \cite{[JD]}.

The second of the graphs that give $3$-simplices
contributes 4 cubes, one for each of the
subforests $\{a,d,e\}$, $\{a,b,d\}$, $\{a,c,d\}$,
and  $\{a,c,e\}$.
The cube corresponding to $\{a,d,e\}$ has a
free minusface given by collapsing the
edge $a$.  So we can disregard this cube.
The remaining 3 cubes join together to
form a solid $3$-ball as described in
the proof of Proposition $10.3$ of \cite{[JD]}, and so can also be
collapsed away. \END \end{pf*}


\begin{thebibliography}{1}

\bibitem{[A-M]}
A. Adem and R. J. Milgram,
{\em Cohomology of Finite Groups},
Springer-Verlag Berlin, Heidelberg 1994.

\bibitem{[B-T]}
G, Baumslag and T. Taylor,
{\em The center of groups with one defining relator},
Math. Ann. 175 (1968) 315-319.

\bibitem{[Br]}
T. Brady,
{\em The integral cohomology of $Out_+(F_3)$},
J. Pure Appl. Algebra 87 (1993) 123-167.

\bibitem{[B]}
K. Brown,
{\em Cohomology of Groups},
Springer-Verlag Berlin, Heidelberg 1982.

\bibitem{[Ch]}
Y. Q. Chen,
{\em Farrell cohomology of automorphism groups
of free groups of finite rank},
Ohio State University Ph.D. dissertation, Columbus, Ohio 1998.

\bibitem{[C]}
M. Culler,
{\em Finite groups of outer automorphisms of a free group},
Contemporary Math 33 (1984) 197-207.

\bibitem{[C-V]}
M. Culler and K. Vogtmann,
{\em Moduli of graphs and automorphisms of free groups},
Invent. Math. 84 (1986) 91-119.

\bibitem{[F]}
F. Thomas Farrell,
{\em An extension of Tate cohomology to a class of infinite groups},
J. Pure Appl. Algebra 10 (1977) 153-161.

\bibitem{[G-M]}
H. H. Glover and G. Mislin,
{\em On the $p$-primary cohomology of $Out(F_n)$ in the $p$-
rank one case},
to appear in J. Pure Appl. Algebra 150 (2000), no. 2.

\bibitem{[H]}
Allen Hatcher,
{\em Homological stability for automorphism groups of free groups},
Comment. Math. Helv. 70 (1995) 39-62.

\bibitem{[H-V]}
A. Hatcher and K. Vogtmann,
{\em Cerf theory for graphs},
J. London Math. Soc. (2) 58 (1998), no. 3, 633-655.

\bibitem{[V]}
A. Hatcher and K. Vogtmann,
{\em Rational homology of $Aut(F_n)$},
Math. Res. Lett. 5 (1998), no. 6, 759-780.

\bibitem{[HN]}
Hans-Werner Henn,
{\em Centralizers of elementary abelian $p$-subgroups, the
Borel construction of the singular locus and applications
to the cohomology of discrete groups},
Topology 36 (1997), no. 1, 271-286.

\bibitem{[JD]}
C. A. Jensen,
{\em Cohomology of $Aut(F_n)$},
Cornell University Ph.D. dissertation, Ithaca, New York 1998.

\bibitem{[K]}
Sava Krstic,
{\em Actions of Finite Groups on Graphs and Related Automorphisms
of Free Groups},
J. Algebra 124 (1989) 119-138.

\bibitem{[K-V]}
Sava Krstic and Karen Vogtmann,
{\em Equivariant outer space and automorphisms of
free-by-finite groups},
Comment. Math. Helv. 68 (1993) 216-262.

\bibitem{[M]}
H. Minkowski,
{\em Zur Theorie der positiven quadratischen Formen},
Crelles J. 101 (1887) 196-202.

\bibitem{[S-V1]}
J. Smillie and K. Vogtmann,
{\em A generating function for the Euler characteristic
of $Out(F_n)$},
J. Pure Appl. Algebra 44 (1987) 329-348.

\bibitem{[S-V2]}
J. Smillie and K. Vogtmann,
{\em Automorphisms of graphs, $p$-subgroups of $Out(F_n)$ and the
Euler characteristic of $Out(F_n)$},
J. Pure Appl. Algebra 49 (1987) 187-200.

\bibitem{[S]}
R. G. Swan,
{\em The $p$-period of a finite group},
Ill. J. Math 4 (1960), 341-346.

\bibitem{[Z]}
B. Zimmerman,
{\em \"Uber Hom\"oomorphismen $n$-dimensionaler Henkelk\"orper und
endliche Erweiterungen von Schottky-Gruppen},
Comment. Math. Helv. 56 (1981) 424-486.

\end{thebibliography}
\end{document}